\documentclass[letterpaper]{amsart}

\usepackage{graphicx}
\usepackage{amssymb, comment, color}
\usepackage{mathrsfs}
\usepackage[active]{srcltx}
\usepackage[colorlinks=true]{hyperref}

\newcommand{\N}{\mathbb{N}}                     
\newcommand{\Z}{\mathbb{Z}}                     
\newcommand{\R}{\mathbb{R}}                     
\newcommand{\C}{\mathbb{C}}                     
\newcommand{\D}{\mathbb{D}}
\renewcommand{\P}{\mathcal{P}}
\newcommand{\Pp}{\mathscr{P}}

\newcommand{\CZ}{{\rm CZ}}

\newtheorem{thm}{\sc Theorem}[section]          
\newtheorem*{thm*}{\sc Theorem}               	
\newtheorem{cor}[thm]{\sc Corollary}        	
\newtheorem*{cor*}{\sc Corollary}        		
\newtheorem{defn}[thm]{\sc Definition}      
\newtheorem{rem}[thm]{\sc Remark}           
\newtheorem{ques}{\sc Question}

\title[Global surfaces of section for Reeb flows]{Global surfaces of section for Reeb flows in dimension three and beyond}

\author{Umberto L. Hryniewicz}
\author{Pedro A. S. Salom\~ao}

\begin{document}

\begin{abstract}
We survey some recent developments in the quest for global surfaces of section for Reeb flows in dimension three using methods from Symplectic Topology. We focus on applications to geometry, including existence of closed geodesics and sharp systolic inequalities. Applications to topology and celestial mechanics are also presented.
\end{abstract}

\maketitle

\begin{center}{\it Dedicated to the memory of Professor Kris Wysocki}\end{center}

\tableofcontents


\section{Introduction}

The idea of a global surface of section goes back to Poincar\'e and the planar circular restricted three-body problem.

\begin{defn}
Let $\phi^t$ be a smooth flow on a smooth closed $3$-manifold $M$. A compact embedded surface $\Sigma \hookrightarrow M$ is a global surface of section for $\phi^t$ if:
\begin{itemize}
\item[(i)] Each component of $\partial\Sigma$ is a periodic orbit of $\phi^t$.
\item[(ii)] $\phi^t$ is transverse to $\Sigma\setminus\partial\Sigma$.
\item[(iii)] For every $p\in M\setminus\partial\Sigma$ there exist $t_+>0$ and $t_-<0$ such that $\phi^{t_+}(p)$ and $\phi^{t_-}(p)$ belong to $\Sigma\setminus\partial\Sigma$.
\end{itemize}
\end{defn}

Any $p\in\Sigma\setminus\partial\Sigma$ has a first return time $\tau(p) = \inf \{ t>0 \mid \phi^t(p)\in \Sigma \}$ and the global dynamical behavior of the flow is encoded in the first return map
\begin{equation}
\psi : \Sigma\setminus\partial\Sigma \to \Sigma\setminus\partial\Sigma \qquad \psi(p) = \phi^{\tau(p)}(p)
\end{equation}
In~\cite{Po} Poincar\'e described annulus-like global surfaces of section for the planar circular restricted three-body problem~(PCR3BP) for certain values of the Jacobi constant and mass ratio. Poincar\'e's global sections motivated his celebrated {\it last geometric theorem}. The associated first return map preserves an area form, extends up to boundary, and satisfies a twist condition in the range of parameters considered. The exciting discovery made by Poincar\'e was that the twist condition implies the existence of infinitely many periodic points, i.e., infinitely many periodic orbits for the PCR3BP. In one stroke Poincar\'e gave a strong push towards a qualitative point of view for studying differential equations, and stated a fixed point theorem intimately connected to the Arnold conjectures and the foundations of Floer Theory.

The recent success of Floer theory and other methods from Symplectic Geometry prompted Hofer to coin the term {\it Symplectic Dynamics}~\cite{BH}. In this note we are concerned with the success of these methods to study Reeb flows in dimension three, with an eye towards applications to geometry.

The presence of a global surface of section is a ``dynamicist paradise''  according to Ghys~\cite{ghys}. It opens the possibility of using methods in two-dimensional dynamics to understand three-dimensional flows. Our first goal is to discuss existence results for global sections. This will be done in section~\ref{sec_existence_results}. After stating Birkhoff's theorem and some results from Schwartzman-Fried-Sullivan theory~\cite{fried,schwartzman,sullivan}, we shall focus on applications of Hofer-Wysocki-Zehnder's theory of pseudo-holomorphic curves~\cite{93,props1,props2,props3}. Particular emphasis is given to Hamiltonian dynamics on convex energy levels. We survey some published (and also some unpublished) results without giving proofs.

Section~\ref{sec_systolic} is devoted to applications to systolic geometry. We will explain how global surfaces of section open the door for Symplectic Dynamics to be used to obtain some sharp systolic inequalities. After a quick review of the rudiments of systolic geometry, we focus on Riemannian two-spheres and on a special case of a conjecture of Viterbo.

In section~\ref{sec_pcr3bp} we present the planar circular restricted three-body problem in more detail. A conjecture due to Birkhoff on the existence of disk-like global surfaces of section for retrograde orbits is discussed.

We intend to convince the reader that there are many positive results for global sections in large classes of flows. However, there are situations where it might be hard to decide whether they exist or not. In sections~\ref{sec_foliations} and~\ref{sec_PB} we discuss results to handle some of these situations for tight Reeb flows on $S^3$. In section~\ref{sec_foliations} we present some of the deep results of Hofer, Wysocki and Zehnder~\cite{fols} concerning the existence of transverse foliations, and its use in the study of Hamiltonian dynamics near critical levels. In section~\ref{sec_PB} we present a Poincar\'e-Birkhoff theorem for tight Reeb flows on $S^3$ proved in~\cite{HMS}. It concerns Reeb flows with a pair of closed orbits exactly as those in the boundary of Poincar\'e's annulus, i.e. forming a Hopf link.

The appendix presents a new proof of the existence of infinitely many closed geodesics on any Riemannian two-sphere, which is alternative to the classical arguments of Bangert~\cite{bangert} and Franks~\cite{franks}. It uses the work of Hingston~\cite{hingston}.

\section{Existence results for global surfaces of section}\label{sec_existence_results}

\subsection{Birkhoff annuli}

Poincar\'e constructed his annulus map for a specific family of systems which are close to integrable\footnote{Angular momentum is preserved in the rotating Kepler problem.}. One of the first statements for a large family of systems which can be quite far from integrable is the following theorem of Birkhoff.

\begin{thm}[Birkhoff~\cite{birkhoff_book}]\label{thm_birkhoff}
Let $\gamma$ be a simple closed geodesic of a positively curved Riemannian two-sphere. Consider the set $A_\gamma$ of unit vectors along $\gamma$ pointing towards one of the hemispheres determined by $\gamma$. Then $A_\gamma$ is a global surface of section for the geodesic flow.
\end{thm}

In other words, every geodesic ray not contained in $\gamma$ visits both hemispheres infinitely often. We call the embedded annulus $A_\gamma$ the {\it Birkhoff annulus}. The family of geodesic flows on positively curved two-spheres is large, making the above statement quite useful. Still, the proof heavily relies on Riemannian geometry and sheds little light on the general problem of finding global sections.

\subsection{Existence results from Schwartzman-Fried-Sullivan theory}

A very general theory to attack the existence problem of global surfaces of section exists, and nowadays goes by the name of Schwartzman-Fried-Sullivan theory; see~\cite{ghys}. It can be powerful since it handles general flows. The drawback is that it provides theorems with hypotheses that are hard to check. This should not be a surprise because the set of all flows on a $3$-manifold is too wild. We pause to review some basic notions. Our discussion is influenced by that of~\cite{ghys} and the work of Fried~\cite{fried}.

Let $\phi^t$ be a smooth flow on a smooth oriented closed $3$-manifold $M$, with infinitesimal generator $X = \partial_t\phi^t|_{t=0}$. Let $L$ be a link consisting of periodic orbits, which we orient by the flow, or the empty set. The major actors here are the $\phi$-invariant Borel probability measures on $M\setminus L$. The set of these measures will be denoted by $\Pp_\phi(M\setminus L)$.

Given $\mu \in \Pp_\phi(M\setminus L)$ and $y \in H^1(M\setminus L;\R)$, we would like to define a real-valued ``intersection number'' $\mu \cdot y \in \R$ as follows.

Firstly, we need to represent $y$ by special closed $1$-forms $\beta^y$ on $M\setminus L$. For every component $\gamma\subset L$ choose a point $q\in\gamma$ and a neighborhood $N_\gamma$ of $\gamma$ equipped with an orientation preserving diffeomorphism $\Psi: N_\gamma \to \R/T\Z \times \D$ such that $\Psi(\phi^t(q)) = (t,0)$. Here $T>0$ is the primitive positive period of $\gamma$. On $N_\gamma\setminus\gamma$ there are coordinates $(t,r,\theta)$ given by $\Psi^{-1}(t,re^{i\theta}) \simeq (t,r,\theta)$, referred to as {\it tubular polar coordinates around~$\gamma$}. There always exists $\beta^y$ such that $\beta^y = p dt+ qd\theta$ on $N_\gamma$ with constants $p,q\in\R$ determined by $y$, for every $\gamma$. Now we are in position to define
\begin{equation}\label{int_number_defn}
\mu \cdot y = \int_{M\setminus L} \beta^y(X) \ d\mu
\end{equation}
The function $\beta^y(X)$ turns out to be bounded, and the integral~\eqref{int_number_defn} depends only on $y$ and $\mu$.

\begin{rem}\label{rem_periodic_orbit_measure}
A $\tau$-periodic orbit $c:\R/\tau \Z \to M\setminus L$ defines $\mu_c \in \Pp_\phi(M\setminus L)$ as the push-forward of normalized Lebesgue measure on $\R/\tau \Z$. If $f\in H_2(M,L;\Z)$ then there is a dual class $y^f\in H^1(M\setminus L;\R)$ that algebraically counts intersections of classes in $H_1(M\setminus L;\Z)$ with $f$. Then $\mu_c \cdot y^f$ equals to the algebraic intersection number of the class of $c$ with $f$ divided by $\tau$. This justifies our terminology.
\end{rem}

Tubular polar coordinates $(t,r,\theta)$ along a component $\gamma\subset L$ trivialize the normal bundle of $\gamma$. This trivialization defines lifts of arguments $\theta(t)\in\R$ of the linearized flow along $\gamma$. If $y$ is cohomologous to $pdt+qd\theta$ near $\gamma$ then the $y$-rotation number of $\gamma$ can be defined as
\begin{equation}
\rho^y(\gamma) = \frac{T}{2\pi} \left( p + q \lim_{t\to+\infty} \frac{\theta(t)}{t} \right)
\end{equation}
It depends only on $y$ and $\gamma$.

A {\it rational} global surface of section is defined as an immersion $j:\Sigma\to M$ of a compact orientable surface such that $j|_{\Sigma\setminus\partial\Sigma}$ is transverse to $X$ and defines a proper embedding $\Sigma\setminus\partial\Sigma \hookrightarrow M\setminus j(\partial\Sigma)$, $j|_{\partial\Sigma}$ is a covering map $\partial \Sigma \to j(\partial\Sigma)$ onto a collection of periodic orbits $j(\partial\Sigma)$, and every trajectory in $M\setminus j(\partial\Sigma)$ hits $j(\Sigma\setminus \partial \Sigma)$ infinitely often in the future and in the past.

The following statements can be proved via Schwartzman-Fried-Sullivan theory.

\begin{thm}\label{thm_SFS1}
If there exists a class $y\in H^1(M\setminus L;\R)$ such that $\mu\cdot y >0$ for all $\mu \in \Pp_{\phi}(M\setminus L)$, and $\rho^y(\gamma)>0$ for all components $\gamma\subset L$, then the union of some components of $L$ bounds a rational global surface of section for $\phi^t$.
\end{thm}

\begin{thm}\label{thm_SFS2}
Suppose that $L$ is null-homologous and let $f\in H_2(M,L;\Z)$ be induced by a Seifert surface for $L$. If $\mu\cdot y^f >0$ for all $\mu \in \Pp_{\phi}(M\setminus L)$ and $\rho^{y^f}(\gamma)>0$ for all components $\gamma\subset L$, then $L$ binds an open book decomposition whose pages are global surfaces of section for $\phi^t$ representing the class $f$. Conversely, when $X$ is $C^\infty$-generic, these assumptions are necessary for $L$ to bound a global surface of section in class $f$.
\end{thm}

See Remark~\ref{rem_periodic_orbit_measure} for the definition of $y^f$. Both statements do not assume $M$ to be a homology sphere. They can be summarized as saying that linking assumptions on invariant measures ensure the existence of global sections. From a topological perspective, Theorem~\ref{thm_SFS2} explains how dynamics can be used to decide when a knot is fibered.

Proofs follow the pattern described in~\cite{ghys}; for Theorem~\ref{thm_SFS2} one needs a refined version of the Hahn-Banach theorem which is not too difficult to prove. Theorem~\ref{thm_SFS1} might produce a fibration over the circle by closed global surfaces of section and the flow is a suspension, this is certainly the case when $L=\emptyset$.

\subsection{Existence results for Reeb flows from Hofer's theory}

From now on we focus on Reeb flows in three dimensions, in order to compare the above statements to results obtained using pseudo-holomorphic curves in symplectic cobordisms. These were introduced by Hofer~\cite{93}. The theory was developed by Hofer, Wysocki and Zehnder~\cite{props1,props2,props3} and later by many other authors.

Consider $\R^4$ with coordinates $(x_1,y_1,x_2,y_2)$ and its standard symplectic form $\omega_0 = \sum_{j=1}^2 dx_j\wedge dy_j$.  Here is a sample result. \\

{\it A periodic orbit of the Hamiltonian flow on a smooth, compact and strictly convex energy level in $(\R^4,\omega_0)$ bounds a disk-like global surface of section if, and only if, it is unknotted and has self-linking number $-1$.} \\

It is instructive to try to understand what is behind such an {\it unconditional} statement. All hypotheses of Theorem~\ref{thm_SFS2} must be somehow verified, but here it seems that there is nothing to be done for unknotted orbits with self-linking number $-1$.  To reconcile the above statement with Theorem~\ref{thm_SFS2}, and to make more statements about Reeb dynamics, we need to review basic notions.

\subsubsection{Some basic notions in Reeb dynamics}

A contact form $\lambda$ on a $3$-manifold $M$ is a $1$-form such that $\lambda\wedge d\lambda$ defines a volume form. Its Reeb vector field $R_\lambda$ is implicitly defined by
\begin{equation}\label{defn_Reeb_vector_field}
d\lambda(R_\lambda,\cdot)=0 \qquad \lambda(R_\lambda)=1
\end{equation}
The distribution $\xi=\ker\lambda$ is a contact structure, the pair $(M,\xi)$ is a contact manifold. More precisely these are the co-orientable contact manifolds since $\lambda$ orients $TM/\xi$. Here we only work with co-orientable contact structures. Contact manifolds are the main objects of study in contact topology.

Our interest is shifted towards dynamics. By a Reeb flow on $(M,\xi)$ we mean one associated to a contact form $\lambda$ on $M$ such that $\xi=\ker\lambda$. Very simple contact manifolds may carry extremely complicated Reeb dynamics. For instance, the standard contact structure $\xi_0$ on $S^3=\{(x_1,y_1,x_2,y_2)\in\R^4\mid x_1^2+y_1^2+x_2^2+y_2^2=1\}$, defined as the kernel of $$ \lambda_0 = \frac{1}{2}\sum_{j=1}^2 x_j dy_j - y_jdx_j $$ restricted to $S^3$, carries rich Reeb dynamics.

A knot is called transverse if at every point its tangent space is transverse to the contact structure. A transverse knot with a Seifert surface has a self-linking number, which is invariant under transverse isotopies. It is defined as follows: choose a non-vanishing section of the contact structure along the Seifert surface, then use this section to push the knot off from itself, and finally count intersections with the Seifert surface.  The vector bundle $(\xi,d\lambda)$ is symplectic and has a first Chern class $c_1(\xi) \in H^2(M;\Z)$. If $c_1(\xi)$ vanishes on $H_2(M;\Z)$ then the self-linking number does not depend on the Seifert surface. The book~\cite{Ge} by Geiges is a nice reference for these concepts.

We end this discussion with a description of the Conley-Zehnder index taken from~\cite[appendix]{fols}. Let $\gamma$ be a periodic trajectory of the flow $\phi^t$ of the Reeb vector field $R_\lambda$, and let $T>0$ be a period of $\gamma$. Since $(\phi^t)^*\lambda=\lambda$, we get a path of $d\lambda$-symplectic linear maps $d\phi^t:\xi_{\gamma(0)} \to \xi_{\gamma(t)}$. The orbit $\gamma$ is called degenerate in period $T$ if $1$ is an eigenvalue of $d\phi^T:\xi_{\gamma(0)} \to \xi_{\gamma(0)}$, otherwise it is called non-degenerate in period~$T$. The contact form $\lambda$ is called non-degenerate when every periodic trajectory is non-degenerate in every period. When $T$ is the primitive period we may simply call $\gamma$ degenerate or non-degenerate accordingly.

Since $T$ is a period, we get a well-defined map $\gamma:\R/T\Z\to M$ which we still denote by $\gamma$ without fear of ambiguity. Choose a symplectic trivialization $\Phi$ of $\gamma^*\xi$. Then the  linearized flow $d\phi^t:\xi_{\gamma(0)} \to \xi_{\gamma(t)}$ gets represented as a path of symplectic matrices $M:\R\to Sp(2)$ satisfying $M(0)=I$, $M(t+T)=M(t)M(T) \ \forall t$. For every non-zero $u\in\R^2$ we write $M(t)u=(r(t)\cos\theta(t),r(t)\sin\theta(t))$ in polar coordinates, for some continuous lift of argument $\theta:\R\to\R$, and define the rotation function $\Delta_M:\R^2\setminus\{(0,0)\} \to \R$ by
\[
\Delta_M(u) = \frac{\theta(T)-\theta(0)}{2\pi}
\]
The image of $\Delta_M$ is a compact interval of length strictly less than $1/2$. The rotation interval $J_M$ is defined as the image of $\Delta_M$. Consider the following function $\tilde\mu(J)$ defined on closed intervals $J$ of length less than $1/2$. If $\partial J\cap\Z=\emptyset$ then
\[
\tilde\mu(J) = \left\{ \begin{aligned} & 2k \qquad \text{if $k\in J$} \\ & 2k+1 \qquad \text{if $J\subset (k,k+1)$} \end{aligned} \right.
\]
If $\partial J\cap\Z\neq\emptyset$ then set $\tilde\mu(J)$ as $\lim_{\epsilon\to0^+}\tilde\mu(J-\epsilon)$. The Conley-Zehnder index can finally be defined as
\[
\CZ^\Phi(\gamma,T) = \tilde\mu(J_M)
\]
We omit the period when it is taken to be the primitive period. If $c_1(\xi)$ vanishes on spheres and $\gamma:\R/T\Z\to M$ is contractible then we write $\CZ^{\rm disk}$ for the index computed with a trivialization that extends to a capping disk.

The Conley-Zehnder index is an extremely important tool. It is related to Fredholm indices of solutions of many of the elliptic equations from Symplectic Topology.

\subsubsection{Existence results of global surfaces of section for Reeb flows}

The following result is somewhat closer to the spirit of Theorem~\ref{thm_SFS2}.

\begin{thm}[\cite{HLS,HS1}]\label{thm_3sphere_intro}
Consider any Reeb flow on $(S^3,\xi_0)$. A periodic orbit $\gamma$ binds an open book decomposition whose pages are disk-like global surfaces of section if it is unknotted, has self-linking number $-1$, satisfies $\rho^{\rm Seifert}(\gamma)>0$, and is linked to all periodic orbits $\gamma':\R/T\Z\to S^3\setminus\gamma$ satisfying either $\CZ^{\rm disk}(\gamma',T)=2$, or $\CZ^{\rm disk}(\gamma',T)=1$ and $\gamma'$ is degenerate in period $T$. Conversely, if the contact form is $C^\infty$-generic then these assumptions are necessary for $\gamma$ to bound a disk-like global surface of section.
\end{thm}

Here $\rho^{\rm Seifert}$ stands for the rotation number computed with respect to a cohomology class dual to a Seifert surface, as explained in Remark~\ref{rem_periodic_orbit_measure}.  The $C^\infty$-generic condition needed in the converse is implied by non-degeneracy of $\gamma$.

We are in position to compare Theorem~\ref{thm_SFS2} with Theorem~\ref{thm_3sphere_intro}. Theorem~\ref{thm_SFS2} is extremely general since it concerns arbitrary flows on arbitrary oriented closed $3$-manifolds, and links of periodic orbits with no knowledge on their Seifert genus. The price to be paid is that it is difficult to be applied since one needs linking with all invariant measures. On the other hand, Theorem~\ref{thm_3sphere_intro} is about Reeb flows on the $3$-sphere with its standard contact structure, that admit a periodic orbit transversely isotopic to a Hopf fiber. However, this class is large enough to provide interesting applications, and the hypotheses only require linking with a much smaller set of invariant measures, namely those coming from certain periodic orbits.

Why is it that for Reeb flows only a much smaller set of invariant measures needs to be analyzed? One reason is the source of compactness in the analysis. In Schwartzman-Fried-Sullivan theory, foliations transverse to the flow are found by looking for closed differential $1$-forms that evaluate positively on the vector field. Differential $1$-forms can be seen as elements in the dual of the space of $1$-currents, and compactness properties of the space of measures (seen as currents) play a key role in the analysis. All of this goes back to Sullivan's work~\cite{sullivan}. In holomorphic curve theory one tries to work directly with the leaves (curves) of the desired foliation. The relevant generalization of Gromov compactness, called the Symplectic Field Theory (SFT) Compactness Theorem~\cite{sftcomp}, tells us that if the desired transverse foliation can not be obtained as we start deforming a given leaf, then some ``misplaced'' curves must exist. The fundamental results of Hofer~\cite{93} provide unlinked periodic orbits.

We continue to describe applications of Hofer's theory. We need the important notion of dynamical convexity.

\begin{defn}[Hofer, Wysocki and Zehnder]
A contact form $\lambda$ on a $3$-manifold $M$ is dynamically convex if $c_1(\ker\lambda)$ vanishes on spheres and contractible periodic Reeb orbits $\gamma:\R/T\Z\to M$ satisfy $\CZ^{\rm disk}(\gamma,T)\geq3$.
\end{defn}

A Reeb flow will be called dynamically convex when it is induced by a dynamically convex contact form. If a periodic orbit $\gamma$ of a Reeb flow on $(S^3,\xi_0)$ is unknotted and has self-linking number $-1$ then $\rho^{\rm Seifert}(\gamma)>0$ turns out to be equivalent to $\CZ^{\rm disk}(\gamma)\geq 3$. The proof of the converse in Theorem~\ref{thm_3sphere_intro} needs genericity assumptions only to guarantee that $\rho^{\rm Seifert}(\gamma)>0$ is necessary. Consequently we get the following particular case.

\begin{thm}[\cite{hryn,hryn2}]\label{thm_hryn_dyn_convex}
Consider a dynamically convex Reeb flow on $(S^3,\xi_0)$. A~periodic orbit $\gamma$ bounds a disk-like global surface of section if, and only if, it is unknotted and has self-linking number $-1$. Moreover, such an orbit binds an open book decomposition whose pages are disk-like global surfaces of section.
\end{thm}

These statements are powered by the following non-trivial input.

\begin{thm}[Hofer, Wysocki and Zehnder~\cite{unknotted}]\label{thm_unknotted_HWZ}
Any Reeb flow on $(S^3,\xi_0)$ has an unknotted periodic orbit with self-linking number $-1$.
\end{thm}

Putting together Theorem~\ref{thm_hryn_dyn_convex} and Theorem~\ref{thm_unknotted_HWZ} we obtain the following remarkable result, which we see as one of the pinnacles of Symplectic Dynamics.

\begin{thm}[Hofer, Wysocki and Zehnder~\cite{convex}]\label{thm_dyn_convex_HWZ}
Any dynamically convex Reeb flow on $(S^3,\xi_0)$ admits a disk-like global surface of section.
\end{thm}

The order in which results were presented is not chronological. Theorem~\ref{thm_dyn_convex_HWZ} was the first, guiding major result in this subject. All our theorems are proved by pushing the methods from~\cite{convex}.

Reeb flows on $(S^3,\xi_0)$ come from smooth, compact hypersurfaces $S$ in $(\R^4,\omega_0)$ that are star-shaped with respect to the origin. This means that $S$ does not go through the origin and every ray from the origin hits it once and transversely. Then $\lambda_0$ restricts to a contact form on $S$ whose Reeb flow parametrizes the characteristic foliation given by the integral lines of $(TS)^{\omega_0}$. In other words, the Reeb flow reparametrizes the Hamiltonian flow on $S$ for any Hamiltonian realizing $S$ as a regular energy level. Moreover, such a Reeb flow on $S$ is smoothly conjugated to a Reeb flow on $(S^3,\xi_0)$. Conversely, every Reeb flow on $(S^3,\xi_0)$ is smoothly conjugated to the Reeb flow of $\lambda_0$ restricted to some $S$.

The following theorem explains the role of convexity, and why the statement from the beginning of this section follows from Theorem~\ref{thm_hryn_dyn_convex}.

\begin{thm}[Hofer, Wysocki and Zehnder~\cite{convex}]\label{thm_convexity_HWZ}
Every Hamiltonian flow on a smooth, compact and strictly convex energy level in $(\R^4,\omega_0)$ is equivalent to a dynamically convex Reeb flow on $(S^3,\xi_0)$.
\end{thm}

We mentioned in the beginning of this introduction that global sections open the door for tools in two-dimensional dynamics. Here is a strong application taken from~\cite{convex}. The return map of the disk obtained from Theorem~\ref{thm_dyn_convex_HWZ} preserves an area form with finite total area. Brouwer's translation theorem provides a periodic orbit simply linked to the boundary of the disk. A strong result due to John Franks~\cite{franks} implies the following statement.

\begin{cor}[Hofer, Wysocki and Zehnder~\cite{convex}]
Hamiltonian flows on smooth, compact and strictly convex energy levels in $(\R^4,\omega_0)$ admit either two or infinitely many periodic orbits.
\end{cor}

The simplest example of convex energy level comes from a pair of uncoupled harmonic oscillators. The two periodic orbits obtained by freezing one oscillator organize the dynamics: they form a Hopf link and both span foliations by disk-like global surfaces of section. In other words, the flow is monotonically rotating around these two orbits. Of course this is a simple integrable system, explicit constructions can be made using an integral. Putting together theorems~\ref{thm_hryn_dyn_convex},~\ref{thm_dyn_convex_HWZ} and~\ref{thm_convexity_HWZ} with Brouwer translation theorem, we easily prove the following corollary. It states that the same holds for every convex energy level, and provides new global surfaces of section that are invisible\footnote{In Theorem~\ref{thm_dyn_convex_HWZ} the boundary of the obtained disk-like global surface of section has Conley-Zehnder index equal to $3$. A linked orbit coming from a fixed point might have arbitrarily large Conley-Zehnder index.} to Theorem~\ref{thm_dyn_convex_HWZ}.

\begin{cor}[\cite{hryn2,HS1}]
Consider the Hamiltonian flow on a smooth, compact and strictly convex energy level in $(\R^4,\omega_0)$. There exists a pair of periodic orbits forming a Hopf link, both of which bind open book decompositions whose pages are disk-like global surfaces of section.
\end{cor}

Symmetries often play a key role in dynamics. Dynamically convex energy levels in $(\R^4,\omega_0)$ which are invariant with respect to the anti-symplectic conjugation $$ \widetilde\rho:(x_1,y_1,x_2,y_2)\mapsto (x_1,-y_1,x_2,-y_2) $$ appear in Celestial Mechanics. The following existence result becomes relevant in this context.

\begin{thm}[Frauenfelder and Kang~\cite{FK}]\label{thm_FK}
If the contact form induced by $\lambda_0$ on a smooth, compact, star-shaped and $\widetilde\rho$-invariant hypersurface in $(\R^4,\omega_0)$ is dynamically convex, then its Reeb flow admits a $\widetilde\rho$-invariant disk-like global surface of section.
\end{thm}

As a corollary, the boundary of the disk is a $\widetilde\rho$-invariant unknotted periodic orbit with self-linking number $-1$. Results of Kang~\cite{Kang} ensure that there are either two or infinitely many $\widetilde\rho$-invariant periodic orbits.

We have addressed Reeb flows on the $3$-sphere with its standard contact structure. How about more general Reeb flows? How close can we get to a contact analogue of Theorem~\ref{thm_SFS2}? Can we recover and generalize Birkhoff's Theorem~\ref{thm_birkhoff}?

A nice feature of Theorem~\ref{thm_SFS2} is that one does not assume $L$ to be fibered, but then one uses the hypotheses to conclude that $L$ is fibered. From a topological perspective this is a strong conclusion. However, from a dynamical perspective the theorem would still be very strong even if we knew that $L$ is fibered because it provides open books whose pages are global surfaces of section! We keep this in mind when looking for a reasonable analogue for Reeb flows.

The notion of fibered link has a contact topological analogue. If $\lambda$ is a contact form and $L$ is a transverse link then the right notion of fibered is that $L$ binds an open book decomposition satisfying
\begin{itemize}
\item[(i)] $d\lambda$ is an area form on each page, and
\item[(ii)] the boundary orientation induced on $L$ by the pages oriented by $d\lambda$ coincides with the orientation induced on $L$ by $\lambda$.
\end{itemize}
Such an open book is said to support the contact structure $\xi = \ker\lambda$. We may call them {\it Giroux open books} because of their fundamental role in the classification of contact structures due to Giroux~\cite{gir02}.

Two more pieces of terminology before our next statement. An open book decomposition is said to be {\it planar} if pages have no genus. A global surface of section will be called {\it positive} if the orientation induced on it by the flow and the ambient orientation turns out to orient its boundary along the flow.

\begin{thm}\label{main1}
Let $(M,\xi)$ be a closed, connected contact $3$-manifold. Let the link $L\subset M$ bind a planar Giroux open book decomposition $\Theta$ of $M$. Denote by $f\in H_2(M,L;\Z)$ the class of a page of $\Theta$, and by $\gamma_1,\dots,\gamma_n$ the components of $L$. Let the contact form $\lambda$ define $\xi$ and realize $L$ as periodic Reeb orbits. Consider the following assertions:
\begin{itemize}
\item[(i)] $L$ bounds a positive genus zero global surface of section for the $\lambda$-Reeb flow representing the class~$f$.
\item[(ii)] $L$ binds a planar Giroux open book whose pages are global surfaces of section for the $\lambda$-Reeb flow and represent the class~$f$.
\item[(iii)] The following hold:
\begin{itemize}
\item[(a)] $\rho^{y^f}(\gamma_k)>0$ for all $k$, where $y^f\in H^1(M\setminus L;\R)$ is dual to $f$.
\item[(b)] Every periodic $\lambda$-Reeb orbit in $M\setminus L$ has non-zero intersection number with $f$.
\end{itemize}
\end{itemize}
Then (iii) $\Rightarrow$ (ii) $\Rightarrow$ (i). Moreover, (i) $\Rightarrow$ (iii) provided the following $C^\infty$-generic condition holds: for every $k$, $\rho^{y^f}(\gamma_k)=0 \Rightarrow \ \text{$\gamma_k$ is hyperbolic}$.
\end{thm}

Theorem~\ref{main1} is fruit of joint work with Kris Wysocki and will be proved in~\cite{HSW}. Our collaboration started in Oberwolfach 2011. It heavily relies on Siefring's intersection theory for holomorphic curves in $4$-dimensional symplectic cobordisms~\cite{siefring}. As explained before, there is no ``mystery'' about the contact topology of $L$. The power of the statement comes from the fact that we do not need linking assumptions about all invariant measures in $M\setminus L$, we only need them for those induced by periodic orbits.

As a first test, note that Birkhoff's Theorem~\ref{thm_birkhoff} follows as a consequence. Indeed, the unit sphere bundle of $S^2$ has a contact form induced by pulling back the tautological $1$-from on $T^*S^2$ via Legendre transform. Reeb flow is geodesic flow. A simple closed geodesic lifts to two closed Reeb orbits, which form a link that binds a supporting open book. Pages are annuli that are isotopic to the Birkhoff annulus. Positivity of the curvature implies that (iii-a) holds. Positivity of the curvature and the Gauss-Bonnet theorem imply that (iii-b) holds. Birkhoff's theorem follows.

Theorem~\ref{main1} has applications to Celestial Mechanics. The following statement is the abstract result needed for these applications. The standard primitive $\lambda_0$ of $\omega_0$ is symmetric by the antipodal map. Identifying antipodal points we obtain $\R P^3 = S^3/\{\pm1\}$. The restriction of $\lambda_0$ to $S^3$ descends to a contact form on $\R P^3$ defining its standard contact structure, still denoted $\xi_0$. The Hopf link $$ \widetilde l_0 = \{(x_1,y_1,x_2,y_2)\in S^3 \mid x_1=y_1=0 \ \text{or} \ x_2=y_2=0\} $$ is antipodal symmetric and descends to a transverse link $l_0$ on $\R P^3$. Any transverse link in $(\R P^3,\xi_0)$ transversely isotopic to $l_0$ will be called a Hopf link. Any transverse knot in $(\R P^3,\xi_0)$ transversely isotopic to a component of $l_0$ will be called a Hopf fiber.

\begin{thm}[\cite{HS2}, \cite{HSW}]\label{thm_annulus_dyn_convex}
Consider an arbitrary dynamically convex Reeb flow on $(\R P^3,\xi_0)$. Any periodic orbit which is a Hopf fiber binds an open book decomposition whose pages are rational disk-like global surfaces of section. Any pair of periodic orbits forming a Hopf link binds an open book decomposition whose pages are annulus-like global surfaces of section.
\end{thm}

In the next section we will explain how the above statement can be applied to Celestial Mechanics. In particular it generalizes Conley's result from~\cite{conley_cpam}.

These techniques have applications to existence of elliptic periodic orbits. A periodic orbit is elliptic if all Floquet multipliers lie in the unit circle. Conjecturally every compact regular and strictly convex energy level in $(\R^{4},\omega_0)$ has an elliptic periodic orbit. Ekeland~\cite{ekeland_elliptic} proved this conjecture in special cases. If the convex level is antipodal symmetric then it has been proved by Dell'Antonio, D'Onofrio and Ekeland~\cite{AOE}. There are results in this direction by Y. Long and his collaborators~\cite{HLongW,LongZhu}, and by Abreu and Macarini~\cite{AM}.

\begin{thm}[\cite{HS2}]\label{thm_elliptic}
Any Reeb flow on $(\R P^3,\xi_0)$ which is sufficiently $C^\infty$-close to a dynamically convex Reeb flow admits an elliptic periodic orbit. This orbit binds a rational open book decomposition whose pages are disk-like global surfaces of section. Its double cover has Conley-Zehnder index equal to~$3$.
\end{thm}

When combined with a result of Harris and Paternain~\cite{HP} relating pinched flag curvatures to dynamical convexity, Theorem~\ref{thm_elliptic} refines the main result of Rademacher from~\cite{Rad}.

\begin{cor}
Consider a Finsler metric on the two-sphere with reversibility~$r$. If all flag curvatures lie in $(r^2/(r+1)^2,1]$ then there exists an elliptic closed geodesic. Moreover, its velocity vector defines a periodic orbit of the geodesic flow that bounds a rational disk-like global surface of section. A fixed point of the return map gives a second closed geodesic.
\end{cor}

The quest for elliptic closed geodesics on a Riemannian two-sphere is an extremely hard problem that goes back to Poincar\'e. In the 1980's Ballmann, Thorbergsson and Ziller~\cite{BTZ} obtained elliptic closed geodesics on a Riemannian $S^n$ under pinching conditions on the curvature. Contreras and Oliveira~\cite{CO} show that any Riemmanian metric on $S^2$ can be $C^2$-approximated by one with an elliptic closed geodesic. Their proof uses a deep result due to Hofer, Wysocki and Zehnder~\cite{fols} asserting the existence of foliations transverse to a non-degenerate Reeb vector field on $(S^3,\xi_0)$, see section~\ref{sec_foliations} for a discussion.

\subsubsection{Dynamical characterization of lens spaces}

We end this section with a topological application. We look for characterizations of contact $3$-manifolds in terms of Reeb dynamics, motivated by early fundamental results of Hofer, Wysocki and Zehnder~\cite{char1,char2}.

Identify $\R^4\simeq\C^2$ by $(x_1,y_1,x_2,y_2) \simeq (z_1=x_1+iy_1,z_2=x_2+iy_2)$ and fix relatively prime integers $p\geq q\geq 1$. In these coordinates we get
\[
\lambda_0 = \frac{i}{4} \sum_{j=1}^2 z_jd\bar z_j - \bar z_jdz_j \qquad \omega_0 = \frac{i}{2} \sum_{j=1}^2 dz_j\wedge d\bar z_j
\]
from where it follows that $\lambda_0$ is invariant under the $\Z_p$ action generated by
\[
(z_1,z_2) \mapsto (e^{i2\pi/p}z_1,e^{i2\pi q/p}z_2)
\]
This action is free on  $S^3 = \{(z_1,z_2)\in \C^2\mid z_1\bar z_1+z_2\bar z_2=1\}$. The lens space $L(p,q)$ is defined as
\[
L(p,q) = S^3/\Z_p
\]
The $3$-sphere is included as the case $p=q=1$, and $\R P^3=L(2,1)$. The contact form induced by $\lambda_0$ on $S^3$ descends to $L(p,q)$ and defines there its so-called standard contact structure, again denoted by $\xi_0$.

A knot $K$ on a closed $3$-manifold $M$ is $p$-unknotted if there is an immersion $u:\D\to M$ such that $u|_{\D\setminus\partial\D}$ defines a proper embedding $\D\setminus\partial\D \to M\setminus K$, and $u|_{\partial\D}$ defines a $p$-covering map $\partial\D \to K$. The map $u$ is called a $p$-disk for $K$. The Hopf fiber $S^1\times 0 \subset S^3$ is $\Z_p$-invariant and descends to the simplest example of a $p$-unknotted knot in $L(p,q)$. The case $p=2$ has the following geometric meaning: if we identify $L(2,1)$ with the unit tangent bundle of the round two-sphere then the velocity vector of a great circle is $2$-unknotted. In the presence of a contact structure a transverse $p$-unknotted knot has a rational self linking number. In the examples given above the knots are transverse and their rational self-linking numbers are equal to $-1/p$. These notions play a role in the following dynamical characterization of standard lens spaces $(L(p,q),\xi_0)$.

\begin{thm}[\cite{char1,char2,HLS}]\label{thm_lens}
Let $(M,\xi)$ be a closed connected contact $3$-manifold, and let $p\geq1$ be an integer. Then $(M,\xi)$ is contactomorphic to some $(L(p,q),\xi_0)$ if, and only if, it carries a dynamically convex Reeb flow with a $p$-unknotted self-linking number $-1/p$ periodic orbit.
\end{thm}

This is a special case of more general statements where linking hypotheses with certain periodic orbits are assumed. The mere existence of a $p$-unknotted self-linking number $-1/p$ periodic orbit implies that $(M,\xi) = (L(p,q),\xi_0)\#(M',\xi')$ for some contact $3$-manifold $(M',\xi')$. The dynamical convexity assumption forces $(M',\xi')=(S^3,\xi_0)$.

Using that $(L(2,1),\xi_0)$ is contactomorphic to the unit sphere bundle of any Finsler metric on $S^2$ we get a geometric application. Consider the set $\mathcal{I}$ of immersions $S^1\to S^2$ with no positive self-tangencies. Two immersions are declared equivalent if they are homotopic through immersions in $\mathcal{I}$. This defines an equivalence relation $\sim$ and an element of $\mathcal{I}/\sim$ will be called a {\it weak flat knot type}. This notion is related to Arnold's $J^+$-theory of plane curves. Note that a closed geodesic of some Finsler metric on $S^2$ has its weak flat knot type well-defined. Let $k_8$ be the weak flat knot type of an eight-like curve, which is a curve with precisely one self-intersection which is transverse. Clearly there are curves representing $k_8$ with an arbitrarily large number of self-intersections.

\begin{thm}[\cite{HS_cambridge}]
If a Finsler two-sphere with reversibility $r$ has flag curvatures in $(r^2/(r+1)^2,1]$ then no closed geodesic represents $k_8$.
\end{thm}

This statement follows from Theorem~\ref{thm_lens}. In fact, the pinching of the curvature forces dynamical convexity (Harris and Paternain~\cite{HP}), and the velocity vector of a closed geodesic of type $k_8$ is unknotted with self-linking number $-1$ in the unit sphere bundle. Since $\R P^3$ is not the $3$-sphere we conclude that such a closed geodesic does not exist.


\section{Global surfaces of section applied to systolic geometry}\label{sec_systolic}

The goal of this section is to explain how Birkhoff's annulus-like global surfaces of section from Theorem~\ref{thm_birkhoff} open the door for symplectic methods to be used to prove sharp systolic inequalities on Riemannian two spheres. And, more generally, how disk-like global surfaces of section can be used in the proof of the perturbative case of Viterbo's conjecture~\cite{viterbo}. The results described here were obtained in collaboration with Alberto Abbondandolo and Barney Bramham~\cite{abhs17,abhs17b,abhs17c,abhs17d}.

\subsection{Origins of systolic geometry}

We start with a quick introduction to systolic geometry. The $1$-systole ${\rm sys}_1(X,g)$ of a closed non-simply connected Riemannian manifold $(X,g)$ is defined as the length of the shortest non-contractible loop. According to Gromov~\cite{gromov_IHES}, {\it ``L\"owner made an amazing discovery around 1949''}, where he refers to the following result.

\begin{thm}[L\"owner]
The inequality $({\rm sys}_1)^2/{\rm Area} \leq 2/\sqrt{3}$ holds for every Riemannian metric on the two-torus. Equality is achieved precisely for the flat torus defined by an hexagonal lattice.
\end{thm}

Another result that can be seen as foundational is Pu's inequality.

\begin{thm}[Pu]
The inequality $({\rm sys}_1)^2/{\rm Area} \leq \pi/2$ holds for every Riemannian metric on $\R P^2$. Equality is achieved precisely for the round geometry.
\end{thm}

Since then the field developed quite a lot. We emphasize Gromov's celebrated paper~\cite{gromov_filling}.

\begin{thm}[Gromov~\cite{gromov_filling}]
There exists a constant $c_n>0$ such that the inequality $({\rm sys}_1)^n/{\rm Volume} \leq c_n$ holds for every closed essential Riemannian manifold of dimension $n$.
\end{thm}

Inequalities such as those in the above theorems will be referred to as {\it systolic inequalities}. Systolic geometry is a huge and active field. One can define higher dimensional systoles, as done by Marcel Berger, and look for systolic and intersystolic inequalities. We found Gromov's lecture notes~\cite{gromov_IHES}, Katz's book~\cite{katz_book} and Guth's papers, such as~\cite{guth1,guth2}, extremely useful.

\subsection{Some problems related to sharp systolic inequalities}

Now we move on to describe some problems. We look for {\it sharp} systolic inequalities. This means that we would like to understand best constants and possibly characterize the geometry in the case of equality, just like in the theorems of L\"owner and Pu.

The case of surfaces is somewhat well studied. However, the two-sphere is left out from the above discussion. To include simply connected manifolds one considers the length $\ell_{\rm min}(X,g)$ of the shortest non-constant closed geodesic of a closed Riemannian manifold $(X,g)$. We introduce the systolic ratio defined by
\begin{equation}
\rho_{\rm sys}(X,g) = \frac{\ell_{\rm min}(X,g)^n}{{\rm Vol}(X,g)} \qquad (n=\dim X)
\end{equation}
and look for the supremum of $\rho_{\rm sys}(S^2,g)$ among all metrics $g$ on $S^2$.

\begin{thm}[Croke~\cite{croke}]\label{thm_croke_upper_bound}
The function $g\mapsto\rho_{\rm sys}(S^2,g)$ is bounded among all Riemannian metrics on $S^2$.
\end{thm}

The best bound seems to be due to Rotman~\cite{rotman}: {\it $\rho_{\rm sys}(S^2,g)\leq 32$ holds for all $g$.} In view of Pu's inequality it is tempting to hope that a round two-sphere $(S^2,g_0)$ maximizes the systolic ratio. Its value is $\rho_{\rm sys}(S^2,g_0)=\pi$. However, the Calabi-Croke sphere shows that
\begin{equation*}
\sup_{(S^2,g)} \rho_{\rm sys}(S^2,g) \geq 2\sqrt{3} > \pi.
\end{equation*}
This is a singular metric constructed by glueing two equilateral triangles along their sides to form a ``flat'' two-sphere. It can be approximated by smooth positively curved metrics with systolic ratio close to $2\sqrt{3}$.

\begin{ques}\label{ques_systolic_two_spheres}
What is the value of $\sup_{(S^2,g)} \rho_{\rm sys}(S^2,g)$? Are there restrictions on the kinds of geometry that approximate this supremum?
\end{ques}

It is conjectured in~\cite{balacheff2} that the answer to Question~\ref{ques_systolic_two_spheres} is $2\sqrt{3}$. In~\cite{balacheff2} Balacheff shows that the Calabi-Croke sphere can be seen as some kind of local maximum if non-smooth metrics with a certain type of singular behavior are included.

A {\it Zoll metric} is one such that all geodesic rays are closed and have the same length. It is interesting that all Zoll metrics on $S^2$ have conjugated geodesic flows, and have systolic ratio equal to $\pi$.

It becomes a natural problem that of understanding the geometry of the function $\rho_{\rm sys}$ near $(S^2,g_0)$. This problem was considered by Babenko and studied by Balacheff. It is tempting to study the function $\rho_{\rm sys}$ with ``variational eyes'' in order to determine its behavior near $(S^2,g_0)$. A serious difficulty is that all Zoll metrics have the same systolic ratio as $(S^2,g_0)$, and Guillemin~\cite{guillemin} found an infinite dimensional family of Zoll metrics through the round geometry. Hence the situation is quite different from the torus or the projective space: one can not hope to determine the metric in terms of systolic ``equalities'' near a local maximizer. In~\cite{balacheff1} Balacheff shows that $(S^2,g_0)$ can be seen as a critical point of $\rho_{\rm sys}$ and conjectured that it is a local maximum. We will refer to this conjecture as the Babenko-Balacheff conjecture.

Contact geometry is a natural set-up to study systolic inequalities. This point of view was advertised and used by \'Alvarez-Paiva and Balacheff~\cite{APB}. See also Michael Hutchings blog post~\cite{Hu_blog}. Let $\alpha$ be a contact form on a closed manifold $M$ of dimension $2n-1$ oriented by $\alpha\wedge(d\alpha)^{n-1}$. We denote by $T_{\rm min}(M,\alpha)$ the minimal period among closed orbits of the Reeb flow. Existence of closed orbits is taken for granted, in dimension three this follows from celebrated result of Taubes~\cite{taubes}. The contact volume of $(M,\alpha)$ is defined as
$$
{\rm Vol}(M,\alpha): = \int_M \alpha \wedge (d \alpha)^{n-1}
$$
and the systolic ratio of $(M,\alpha)$ as
$$
\rho_{\rm sys}(M,\alpha) : = \frac{T_{\rm min}(M,\alpha)^n}{{\rm Vol}(M,\alpha)}
$$
Note that $\rho_{\rm sys}(M,\alpha)$ is invariant under re-scalings of $\alpha$.

To see the connection to systolic geometry, consider a Riemannian $n$-manifold $(X,g)$. Legendre transform pulls-back the tautological $1$-form from $T^*X$ to a $1$-form on $TX$ that restricts to a contact form $\alpha_g$ on the unit sphere bundle $T^1X$. Since the Reeb flow of $\alpha_g$ is the geodesic flow of $g$, we get $T_{\rm min}(T^1X,\alpha_g) = \ell_{\rm min}(X,g)$. It turns out that ${\rm Vol}(X,g)$ and ${\rm Vol}(T^1X,\alpha_g)$ are proportional by a constant depending only on~$n$. Hence we find $C_n>0$ such that
\begin{equation*}
\rho_{\rm sys}(T^1X,\alpha_g) = C_n \ \rho_{\rm sys}(X,g)
\end{equation*}
for every Riemannian metric $g$ on $X$.

A Zoll contact form is one such that all Reeb trajectories are periodic and have the same period. These are usually called regular in the literature, but we prefer the term Zoll in view of the above connection to the Riemannian case.

A convex body in $\R^{2n}$ is a compact convex set with non-empty interior. In~\cite{viterbo} Viterbo conjectured that
\begin{equation}\label{viterbo_full}
\frac{c(K)^n}{n!{\rm Vol}(K)} \leq 1
\end{equation}
holds for every convex body $K\subset \R^{2n}$ and every symplectic capacity $c$, where ${\rm Vol}(K)$ denotes euclidean volume. Viterbo used John's ellipsoid to show that the left-hand side of~\eqref{viterbo_full} is bounded by a constant depending on $n$. In~\cite{AAMO} Artstein-Avidan, Milman and Ostrover showed that the constant can be taken independently of $n$.

In~\cite{AAKO} Artstein-Avidan, Karasev and Ostrover established the connection between~\eqref{viterbo_full} and Mahler conjecture from convex geometry. A centrally symmetric convex body $C\subset \R^n$ has a polar body $C^o=\{y\in\R^n \mid \left<y,x\right> \leq1 \ \forall x\in C\}$. The Mahler conjecture asserts that
\begin{equation}\label{mahler}
{\rm Vol}(C\times C^o) \geq \frac{4^n}{n!}
\end{equation}
for every centrally symmetric convex body $C\subset \R^n$. The main result of~\cite{AAKO} asserts that for every centrally symmetric convex body $C\subset \R^n$ the Hofer-Zehnder capacity of $C\times C^o$ equals $4$. Then~\eqref{mahler} follows from~\eqref{viterbo_full} applied to the Hofer-Zehnder capacity and to convex bodies of the form $C\times C^o$, with a centrally symmetric $C$.

We end by discussing a special case of Viterbo's conjecture. Let $K$ be a convex body in $(\R^{2n},\omega_0)$ with smooth and strictly convex boundary, with the origin in its interior. Denote by $\iota:\partial K \to \R^{2n}$ the inclusion map, and by $\lambda_0$ the standard Liouville form $\lambda_0 = \frac{1}{2} \sum_{j=1}^n x_jdy_j-y_jdx_j$. Then $\iota^*\lambda_0$ is a contact form on $\partial K$. In~\cite{HZ} it is claimed that the Hofer-Zehnder capacity of $K$ is equal to $T_{\rm min}(\partial K,\iota^*\lambda_0)$. In this case~\eqref{viterbo_full} is restated as
\begin{equation}\label{viterbo_particular_case}
\rho_{\rm sys}(\partial K,\iota^*\lambda_0) \leq 1
\end{equation}
which is supposed to be an equality if, and only if, $\iota^*\lambda_0$ is Zoll.

\subsection{The role of global surfaces of section}


Recall that for $\delta\in(0,1]$, a positively curved closed Riemannian manifold is said to be $\delta$-pinched if the minimal and maximal values $K_{\min},K_{\rm max}$ of the sectional curvatures satisfy $K_{\rm min}/K_{\rm max}\geq\delta$. On a positively curved two-sphere we write $\ell_{\rm max}$ for the length of the longest closed geodesic without self-intersections. Note that $\ell_{\rm max}$ is finite.

\begin{thm}[\cite{abhs17}]\label{thm_abhs_two-sphere}
If $(S^2,g)$ is $\delta$-pinched for some $\delta > (4+\sqrt{7})/8 = 0.8307...$ then
\begin{equation}
\ell_{\rm min}(S^2,g)^2 \leq  \pi {\rm Area}(S^2,g) \leq \ell_{\rm max}(S^2,g)^2
\end{equation}
Moreover, any of these inequalities is an equality if, and only if, the metric is Zoll.
\end{thm}

This first inequality confirms the Babenko-Balacheff conjecture on an explicit and somewhat large $C^2$-neighborhood of the round geometry. It seems that the upper bound involving $\ell_{\rm max}$ was not known before.

We discuss some related problems before explaining the role of global surfaces of section in the proof of Theorem~\ref{thm_abhs_two-sphere}. The pinching constant $\delta$ seems to be a helpful parameter. For instance, one could consider the non-increasing bounded (Theorem~\ref{thm_croke_upper_bound}) function
\begin{equation*}
\rho:(0,1] \to \R \qquad \rho(\delta) = \sup \{ \rho_{\rm sys}(S^2,g) \mid (S^2,g) \ \text{is $\delta$-pinched} \}
\end{equation*}
to study the positively curved case.


\begin{ques}
Is it true that $\rho(1/4)=\pi$? What does the graph of $\rho(\delta)$ look like?
\end{ques}

The Calabi-Croke sphere shows that $\lim_{\delta\to0^+} \rho(\delta) \geq 2\sqrt{3}$. Theorem~\ref{thm_abhs_two-sphere} implies that $\rho(\delta) = \pi$ for all $\delta>(4+\sqrt{7})/8$. One must try to understand among which metrics does the round metric maximize systolic ratio. Assuming positive curvature it might be reasonable to expect that $\inf\{\delta \mid \rho(\delta)=\pi\}\leq1/4$.

If curvature assumptions are dropped then the situation might be much harder. What about symmetry assumptions? Here is a result in this direction that answers a question by \'Alvarez-Paiva and Balacheff.

\begin{thm}\label{thm_revolution}
Inequality $\rho_{\rm sys}\leq \pi$ holds for every sphere of revolution, with equality precise when the metric is Zoll.
\end{thm}

%
%
%

Global surfaces of section show up in the proofs of theorems~\ref{thm_abhs_two-sphere} and~\ref{thm_revolution} to connect systolic inequalities to a quantitative fixed point theorem for symplectic maps of the annulus. We outline the proof to make this point precise.

Let $(S^2,g)$ be $\delta$-pinched. If $\delta>1/4$ then $\ell_{\rm min}$ is only realized by simple closed geodesics. Let $\gamma$ be a closed geodesic of length $\ell_{\rm min}$. By Theorem~\ref{thm_birkhoff} the Birkhoff annulus $A_\gamma$ is a global surface of section. Let $\lambda$ be the $1$-form on $A_\gamma$ given by restricting the contact form $\alpha_g$. Then $d\lambda$ is an area form on the interior of $A_\gamma$, and vanishes on $\partial A_\gamma$. The total $d\lambda$-area is $2\ell_{\rm min}$.

The first return map $\psi$ and the first return time $\tau$ are defined on the interior of $A_\gamma$, but it turns out that they extend smoothly to~$A_\gamma$. Moreover, $\psi$ preserves boundary components. Santal\'o's formula reads
\begin{equation}\label{santalo_formula}
2\pi{\rm Area}(S^2,g) = \int_{T^1S^2} \alpha_g\wedge d\alpha_g = \int_{A_\gamma} \tau \ d\lambda
\end{equation}
Since $\psi$ preserves the $2$-form $d\lambda$, it follows that $\psi^*\lambda-\lambda$ is closed.

We now need to consider lifts of $\psi$ to the universal covering of $A_\gamma$. If $\psi$ admits a lift in the kernel of the {\rm FLUX} homomorphism then $\psi^*\lambda-\lambda$ is exact. The unique primitive $\sigma$ of $\psi^*\lambda-\lambda$ satisfying
\begin{equation*}
\sigma(p) = \int_p^{\psi(p)} \lambda \ \ \ \ \forall p\in\partial A_\gamma
\end{equation*}
is called the {\it action of $\psi$}. Here the integral is taken along the boundary according to the lift with zero FLUX. The Calabi invariant is defined as
\begin{equation*}
{\rm CAL}(\psi) = \frac{1}{\int_{A_\gamma}d\lambda} \int_{A_\gamma} \sigma \ d\lambda = \frac{1}{2\ell_{\rm min}} \int_{A_\gamma} \sigma \ d\lambda
\end{equation*}
Of course, we need to worry about whether $\psi$ admits a lift of zero FLUX, but this follows from reversibility of the geodesic flow.

It is a very general fact that $\tau$ is also a primitive of $\psi^*\lambda-\lambda$.  Toponogov's theorem proves that if $\delta>1/4$ then
\begin{equation}\label{returntime_and_action}
\tau = \sigma + \ell_{\rm min}
\end{equation}
Combining~\eqref{returntime_and_action} with~\eqref{santalo_formula} we finally get
\begin{equation}\label{systolic_proof_identity}
2\pi {\rm Area}(S^2,g) = 2(\ell_{\rm min})^2 + 2\ell_{\rm min} {\rm CAL}(\psi)
\end{equation}
Equations~\eqref{returntime_and_action} and~\eqref{systolic_proof_identity} should be seen as some kind of dictionary between geometry and dynamics: action corresponds to length, Calabi invariant corresponds to area.

We are now in position to make the link to the quantitative fixed point theorem and conclude the argument. Roughly speaking, the theorem states: \\

{\it If $\psi$ admits a generating function (of a specific kind), ${\rm CAL}(\psi)\leq0$, and $\psi\neq id$ if ${\rm CAL}(\psi)=0$, then there exists a fixed point $p_0$ satisfying $\sigma(p_0)<0$.} \\

Arguing indirectly, suppose that either $\pi{\rm Area} < (\ell_{\rm min})^2$, or $\pi{\rm Area} = (\ell_{\rm min})^2$ and $g$ is not Zoll. It follows from~\eqref{systolic_proof_identity} and a little more work that either ${\rm CAL}(\psi)<0$, or ${\rm CAL}(\psi)=0$ and $\psi$ is not the identity. Toponogov's theorem comes into play again to show that $\psi$ admits the required generating function provided $\delta>(4+\sqrt{7})/8$. The fixed point theorem applies to give a fixed point of negative action. By~\eqref{returntime_and_action} this fixed point corresponds to a closed geodesic of length strictly smaller than $\ell_{\rm min}$. This contradiction finishes the proof.

The above argument reveals how global surfaces of section can serve as bridge between systolic geometry and symplectic dynamics. The same strategy proves a special case of Viterbo's conjecture in dimension $4$.

\begin{thm}[\cite{abhs17b}]\label{thm_ABHS2}
There exists a $C^3$-neighborhood $\mathcal{U}$ of the space of Zoll contact forms on $S^3$ such that $$ \alpha \in \mathcal{U} \Rightarrow \rho_{\rm sys}(S^3,\alpha) \leq 1 $$ with equality if, only if, $\alpha$ is Zoll.
\end{thm}

The proof again strongly relies on global surfaces of sections. Namely, if a contact form is $C^3$-close to the standard contact form $\lambda_0$ then its Reeb flow admits a disk-like global surface of section whose first return map extends up to the boundary and is $C^1$-close to the identity. We have a dictionary between maps and flows just as in the proof of Theorem~\ref{thm_abhs_two-sphere}: contact volume corresponds to Calabi invariant, return time corresponds to action. The quantitative fixed point theorem applies to give the desired conclusion. 

One could see the constants in sharp systolic inequalities for Riemannian surfaces as invariants. Similarly, one could hope to construct contact invariants from sharp systolic inequalities for contact forms, see~\cite{Hu_blog}. The following statement shows that this is not possible in dimension three: systolic inequalities are not purely contact topological phenomena. For example, inequalities such as~\eqref{viterbo_particular_case} must depend on the convexity assumption.

\begin{thm}[\cite{abhs17b,abhs17c}]\label{thm_ABHS3}
For every co-orientable contact $3$-manifold $(M,\xi)$ and every $c>0$ there exists a contact form $\alpha$ on $M$ satisfying $\xi=\ker\alpha$ and $\rho_{\rm sys}(M,\alpha)> c$.
\end{thm}

Hofer, Wysocki and Zehnder~\cite{char2,convex} introduced the notion of dynamically convex contact forms, see section~\ref{sec_existence_results} for a detailed discussion. It plays a crucial role in the construction of global surfaces of section (theorems~\ref{thm_hryn_dyn_convex},~\ref{thm_dyn_convex_HWZ}). Dynamical convexity is automatically satisfied on the boundary of a smooth convex body if the boundary is strictly convex. It becomes relevant to decide whether~\eqref{viterbo_particular_case} holds for dynamically convex contact forms on $S^3$.

\begin{thm}[\cite{abhs17d}]\label{thm_ABHS4}
Given any $\epsilon>0$ there exists a dynamically convex contact form $\alpha$ on $S^3$ such that $\rho_{\rm sys}(S^3,\alpha) > 2-\epsilon$.
\end{thm}

In \cite{abhs17d}, a narrow connection between high systolic ratios and negativity of Conley-Zehnder indices is quantified.

Observe that Theorem~\ref{thm_ABHS4} implies that either Viterbo's conjecture is not true, or there exists a dynamically convex contact form on $S^3$ whose Reeb flow is not equivalent to the Reeb flow on a strictly convex hypersurface of $(\R^4,\omega_0)$. Unfortunately we can not decide which alternative holds. It also proves that there are smooth compact star-shaped domains $U$ in $(\R^4,\omega_0)$ with the following property: $\partial U$ is dynamically convex and the value $c(U)$ of any capacity realized as the action of some closed characteristic on $\partial U$ is strictly larger than the Gromov width of $U$.

Global surfaces of section continue to play essential role in the proofs of Theorem~\ref{thm_ABHS3} and Theorem~\ref{thm_ABHS4}. Both start by constructing global sections for certain Reeb flows with well-controlled return maps. Then the Reeb flows are modified by carefully changing the return maps in order to make the systolic ratio increase.


We end this section by exploring a bit the theme of the quantitative fixed point theorem. The hypothesis that $\psi$ admits a single generating function is essential, the theorem fails in general. If we note that ${\rm CAL}(\psi)$ is the space average of the action function $\sigma$ then we are led to extrapolate the statement with ``ergodic eyes''. If the ergodic theorem could be applied then the Calabi invariant would be equal to the time average of the action along most orbits. In view of Poincar\'e's recurrence one could hope that ${\rm CAL}(\psi)$ is related to the average action of some periodic point $p$, which is defined as $$ \bar\sigma(p) = \frac{1}{k} \sum_{j=0}^{k-1} \sigma(\psi^j(p)) \qquad \qquad (\text{$k$ is the period of $p$}) $$ This is the content of a difficult result due to Hutchings~\cite{Hu_diskmaps} for disk-maps which we now describe.

Consider the unit disk $\D = \{x+iy \in\C \mid x^2+y^2\leq 1\}$ equipped with the normalized standard area form $\omega = \frac{1}{\pi} dx\wedge dy$, and let $\psi:\D\to\D$ be a diffeomorphism satisfying $\psi^*\omega=\omega$. Assume that there exist numbers $\epsilon>0$ and $\rho\in\R$ such that $r\in(1-\epsilon,1] \Rightarrow \psi(re^{i\theta})=re^{i(\theta+2\pi\rho)}$. Choose a primitive $\lambda$ of $\omega$ satisfying $\lambda(\partial_\theta|_{\partial\D})=1/2\pi$, for instance $\frac{1}{2\pi}(xdy-ydx)$ does the job. If $\sigma:\D\to\R$ is the unique primitive of $\psi^*\lambda-\lambda$ such that $\sigma|_{\partial\D}\equiv\rho$ then set the Calabi invariant of $(\psi,\rho)$ to be its space average
\[
{\rm CAL}(\psi,\rho) = \int_\D \sigma \omega
\]
The dependence on $\rho$ needs to be made explicit because $\psi|_{\partial\D}$ has infinitely many other lifts.

\begin{thm}[Hutchings~\cite{Hu_diskmaps}]
If ${\rm CAL}(\psi,\rho)<\rho$ then for every $\delta>0$ there exists a periodic point $p$ of $\psi$ such that $$ \frac{1}{k} \sum_{j=0}^{k-1} \sigma(\psi^j(p)) < {\rm CAL}(\psi,\rho)+\delta \qquad \text{ where $k$ is a period of $p$.} $$
\end{thm}

The proof uses a global surface of section to embed $\psi$ as a first return map of some Reeb flow on the tight $3$-sphere. The heavy machinery of Embedded Contact Homology and Seiberg-Witten invariants come into play to prove the statement. We hope that the general principle that the Calabi invariant is related to mean actions of periodic orbits extends to more general Hamiltonian systems, this seems to be an interesting topic of research: of course not many relations are expected in general, but they do exist in many special geometric set-ups.

\section{The planar circular restricted three-body problem}\label{sec_pcr3bp}

The three-body problem is that of understanding the motion of three massive particles which attract each other according to Newton's law of gravitation. We follow Moser and Zehnder~\cite{MZ} for basic notation and calculations. The writings of Chenciner, such as~\cite{chenciner}, are recommended for a more comprehensive understanding of the topic.

Denoting positions and masses by $z_1,z_2,z_3$ and $m_1,m_2,m_3$, respectively, the equations of motion become
\begin{align}
& \ddot z_1 = m_2 \frac{z_2-z_1}{|z_2-z_1|^3} + m_3 \frac{z_3-z_1}{|z_3-z_1|^3} \label{eq_12} \\
& \ddot z_2 = m_1 \frac{z_1-z_2}{|z_1-z_2|^3} + m_3 \frac{z_3-z_2}{|z_3-z_2|^3} \label{eq_23} \\
& \ddot z_3 = m_1 \frac{z_1-z_3}{|z_1-z_3|^3} + m_2 \frac{z_2-z_3}{|z_2-z_3|^3} \label{eq_13}
\end{align}
in suitably normalized units. Setting $m_3=0$ in~\eqref{eq_12}-\eqref{eq_23} one obtains the {\it restricted} three-body problem, where the first two particles (primaries) move according to the two-body problem and the third particle (satellite) moves according to~\eqref{eq_13}. Requiring that all particles move on a plane, as we do here, one obtains the {\it planar} problem; from now on the $z_k$ belong to the complex plane $\C$. The relative position $\zeta=z_2-z_1$ solves Kepler's equation $\ddot\zeta=-(m_1+m_2)\zeta/|\zeta|^3$ and the adjective {\it circular} refers to the case when $\zeta$ describes a circle.

In inertial coordinates where the center of mass rests at the origin one gets $z_1=r_1e^{i\omega t}$ and $z_2=-r_2e^{i\omega t}$ for some $\omega$, where $r_1,r_2>0$ satisfy $m_1r_1-m_2r_2=0$ and $(r_1+r_2)^3\omega^2=m_1+m_2$. For definiteness we assume that $\omega>0$. In rotating (non-inertial) coordinates the position $q(t)\in\C$ of the satellite relative to the second primary is given by $z_3(t) = (q(t)-r_2)e^{i\omega t}$, from where it follows that
\begin{equation}\label{equation_for_u}
\ddot q + 2i\omega \dot q - \omega^2(q-r_2) = - m_1\frac{q-r_1-r_2}{|q-r_1-r_2|^3} - m_2\frac{q}{|q|^3}
\end{equation}
From this point on we fix $\omega$ and the total mass $m_1+m_2$, so that $m_1,m_2,r_1,r_2$ are functions of the mass ratio
\begin{equation}
\mu = \frac{m_1}{m_1+m_2}
\end{equation}
As is well known,~\eqref{equation_for_u} can be put in Hamiltonian form. If we set $p = \dot q + i\omega(q-r_2)$ and consider
\begin{equation}\label{hamiltonian_pcr3bp}
H_\mu(q,p) = \frac{1}{2}|p|^2 + \omega\left<q-r_2,ip\right> - \frac{m_1}{|q-r_1-r_2|} - \frac{m_2}{|q|}
\end{equation}
then~\eqref{equation_for_u} becomes Hamilton's equations
\begin{equation}\label{hamilton_eqns}
\dot q = \nabla_pH_\mu \qquad \dot p = -\nabla_qH_\mu
\end{equation}
Of course, this is a special feature of the circular case, in other cases one ends up with a non-autonomous system. After completing squares, one writes
\begin{equation*}
\begin{aligned}
H_\mu(q,p) &= \frac{1}{2}|\omega(q-r_2)+ip|^2 - \frac{1}{2}|\omega(q-r_2)|^2 - \frac{m_1}{|q-r_1-r_2|} - \frac{m_2}{|q|} \\
&= \frac{1}{2}|\omega(q-r_2)+ip|^2 - U_\mu(q)
\end{aligned}
\end{equation*}
where
\begin{equation}
U_\mu(q) = \frac{1}{2}|\omega(q-r_2)|^2 + \frac{m_1}{|q-r_1-r_2|} + \frac{m_2}{|q|}
\end{equation}
is the {\it effective potential}. The function $U_\mu$ has five critical points, the Lagrange points. A sublevel set below its lowest critical value defines three {\it Hill regions}, two of which are bounded while the third is a neighborhood of $\infty$. Each bounded Hill region is topologically a disk and contains a primary, namely, one is a neighborhood of the origin and the other is a neighborhood of $r_1+r_2$. The boundaries of the Hill regions are called {\it ovals of zero velocity}, since there we have $\omega(q-r_2)=-ip \Leftrightarrow \dot q=0$. From now on we restrict to {\it subcritical cases}, i.e. energy levels $H_\mu=-c$ where $-c$ is below the lowest critical value of $H_\mu$. 

Summarizing, we want to study a $1$-parameter family of dynamical systems determined by Hamiltonians $\{H_\mu\}_{\mu\in(0,1)}$ as in~\eqref{hamiltonian_pcr3bp}, on subcritical energy levels. The complete analysis of the entire family is extremely difficult.

Following Poincar\'e, mathematicians first tried to understand the limiting behavior as $\mu\to0^+$ or as $\mu\to1^-$. The limit as $\mu\to0^+$ is in some ways better behaved then the limit $\mu\to1^-$, but sometimes it is just the other way around. In the limit $\mu=0$ the system describes the so-called rotating Kepler problem, where all mass is concentrated at the origin. The boundary of the bounded Hill region about the origin converges to a circle of definite radius. As $\mu\to1^-$ the bounded Hill region about the origin collapses, and we face a somewhat more singular situation. The case $\mu\sim1$ is called the {\it lunar case}, since it is supposed to approximate the motion of a moon near a planet going around the Sun in an almost circular orbit. The main object of this discussion are certain special periodic orbits, which we now recall.

\begin{defn}
A retrograde orbit is a periodic orbit $t\mapsto (q(t),p(t))$ such that $q(t)$ is in the Hill region about the origin, and describes a curve without self-intersections with winding number $-1$ around the origin. Analogously, a direct orbit is a periodic orbit $t\mapsto (q(t),p(t))$ such that $q(t)$ is in the Hill region about the origin, and describes a curve without self-intersections with winding number $+1$ around the origin.
\end{defn}

For $\mu\sim 0$ the existence of retrograde and direct orbits inside the Hill region about the origin follows from a continuation argument starting from circular orbits of the rotating Kepler problem. For this argument to be valid one needs the ratio $\omega/\alpha$ to satisfy a certain genericity condition.

For $\mu\sim 1$ the existence of a retrograde orbit inside the Hill region about the origin follows from Birkhoff shooting argument~\cite{B}, while the existence of a direct orbit is a difficult problem that has eluded mathematicians for more than a century. It is an approximation of the problem of finding the trajectory of a moon around a planet in the solar system which moves in the plane spanned by the sun and the planet and goes around the planet in the same sense as the planet goes around the sun. Our moon behaves in this way.

The difficulty for finding direct orbits led Birkhoff to consider the following strategy in~\cite[section~19]{B}. Firstly one should try to find a disk-like global surface of section bounded by a (doubly covered) retrograde orbit. For this to make sense collision orbits need to be regularized. Secondly, due to preservation of an area form with finite total area, one can apply Brouwer's theorem to the first return map to find a fixed point that should correspond to a direct orbit. Two main difficulties are: (1) for an arbitrary mass ratio it is hard to find global surfaces of section, and (2) it might be hard to check that the fixed point corresponds to a direct orbit. The following is extracted from~\cite[section~19]{B}: \\

{\it ``This state of affairs seems to me to make it probable that the restricted problem of three bodies admit of reduction to the transformation of a discoid into itself as long as there is a closed oval of zero velocity about J, and that in consequence there always exists at least one direct periodic orbit of simple type.''} \\

More recently this has been called a conjecture, which perhaps should be read as following: For any value of $\mu$ and any subcritical energy value, there must be a way of finding a disk-like global surface of section in order to understand the movement of the satellite inside the Hill region about the origin. To implement the strategy of Birkhoff this disk should be spanned by the retrograde orbit, in particular fixed points could be good candidates for direct orbits. Again, all this only makes sense if collision orbits are regularized, and it helps if a regularization procedure is not specified.

Note that the smallest critical value of $H_\mu$ converges to
\begin{equation}
-\frac{3}{2}((m_1+m_2)\omega)^{2/3}
\end{equation}
both when $\mu\to0^+$ or $\mu\to1^-$. Here is a good point to state and discuss our result concerning Birkhoff's conjecture before describing specific regularization procedures.

\begin{thm}\label{thm_birkhoff_conj}
For every $c>\frac{3}{2}((m_1+m_2)\omega)^{2/3}$ there exists $\epsilon>0$ such that the following holds.
\begin{itemize}
\item[(a)] If $1-\mu<\epsilon$ and collisions are regularized via Levi-Civita regularization, then the double cover of every retrograde orbit inside the Hill region about the origin bounds a disk-like global surface of section. Moreover, if we quotient by antipodal symmetry then this disk descends to a rational global disk-like global surface of section.
\item[(b)] If $\mu<\epsilon$ and collisions are regularized via Moser regularization, then the the double cover of every retrograde orbit inside the Hill region about the origin bounds a rational disk-like global surface of section.
\end{itemize}
\end{thm}

A few remarks are in order. Results of Albers, Fish, Frauenfelder, Hofer and van Koert from~\cite{AFFHvK} imply that if $1-\mu$ is small enough then Levi-Civita regularization lifts the Hamiltonian flow on the corresponding component of $H_\mu^{-1}(-c)$ to the characteristic flow on a strictly convex hypersurface $\widetilde\Sigma_{\mu,c}$, up to time reparametrization. Moreover, $\widetilde\Sigma_{\mu,c}$ is antipodal symmetric and each state is represented twice as a pair of antipodal points. Results from~\cite{convex} apply and give disk-like global surfaces of section in $\widetilde\Sigma_{\mu,c}$. Statement (a) above says that there is such a global section in $\widetilde \Sigma_{\mu,c}$ spanned by the lift of every doubly covered retrograde orbit, and that it descends to a global section in the quotient $\Sigma_{\mu,c} = \widetilde\Sigma_{\mu,c}/\{\pm1\}$. In other words, the interior of the disk does not contain pairs of antipodal points, and we end up with a rational disk-like global surface of section in $\Sigma_{\mu,c}\simeq \R P^3$. Note that there is no ambiguity in the representation of points of the original phase space. If $\mu=0$ then Moser regularization applies to the rotating Kepler problem to compactify the Hamiltonian flow on $H_\mu^{-1}(-c)$ to the characteristic flow on a fiberwise starshaped hypersurface $\Sigma_{\mu,c}$ inside $TS^2$, up to time reparametrization. A proof of this statement can be found in the paper~\cite{AFvKP} where the contact-type property of energy levels of the PCR3BP is studied. Again we have $\Sigma_{\mu,c} \simeq \R P^3$. Statement (b) above says that the double cover of every retrograde orbit bounds a rational disk-like global surface of section in $\Sigma_{\mu,c}$. A proof in this case would rely on the dynamical convexity obtained in~\cite{AFFvK} for $\mu=0$.

Fix again $c>\frac{3}{2}((m_1+m_2)\omega)^{2/3}$. As observed above, if $\mu\to0^+$ then we can regularize the component of $H_\mu^{-1}(-c)$ corresponding to the Hill region about the origin via Moser regularization. The result is a Reeb flow on $(\R P^3,\xi_0)$. If $\mu\to1^-$ then we can regularize the component of $H_\mu^{-1}(-c)$ corresponding to the Hill region about the origin via Levi-Civita regularization. The result is a Reeb flow on $(S^3,\xi_0)$ but states are represented twice, by pairs of antipodal. Taking the quotient by the antipodal map the result is again a Reeb flow on $(\R P^3,\xi_0)$. In the first case the limiting system $\mu=0$ is dynamically convex, this is the result of~\cite{AFFvK}. In the second case dynamical convexity is the result of~\cite{AFFHvK}. Hence, for $\mu$ close to $0$ or $1$ we can always apply Theorem~\ref{thm_elliptic} and obtain a pair of periodic orbits which are $2$-unknotted  and have self-linking number $-1/2$. These orbits are transversely isotopic to (a quotient of) a Hopf link. We note that the retrograde and direct orbits obtained from Poincar\'e's continuation method as $\mu\to0^+$ form exactly the same kind of link. Hence Theorem~\ref{thm_elliptic} can be seen as a generalization of this well-known phenomenon in the PCR3BP. Moreover, Theorem~\ref{thm_annulus_dyn_convex} can also be applied and an annulus-like global surface of section is obtained.

We end with a sketch of proof of (a) in Theorem \ref{thm_birkhoff_conj}. It is harmless to assume that $\omega = r_1+r_2=m_1+m_2= 1$ and hence the Hamiltonian in \eqref{hamiltonian_pcr3bp} takes the form
$$
H_\mu(q,p) = \frac{1}{2}|p|^2 + \left<q-\mu,ip\right> - \frac{\mu}{|q-1|} - \frac{1-\mu}{|q|}.
$$ It depends on the mass ratio $\mu=m_1=r_2\in (0,1)$.

The component $\dot \Sigma_{\mu,c}\subset H_\mu^{-1}(-c)$ which projects to the Hill region surrounding $0\in \C$ contains collision orbits. These orbits are regularized with the aid of Levi-Civita coordinates $(v,u)\in \C \times \C$
\begin{equation}\label{change}
q=2v^2 \mbox{  and  } p=-\frac{u}{\bar v},
\end{equation} which are symplectic up to a constant factor.

 The regularized Hamiltonian is
 \begin{equation}\label{Kmc}
 \begin{aligned}
 K_{\mu,c}(v,u) & := |v|^2(H_\mu(p,q)+c)\\ & = \frac{1}{2}|u|^2+2|v|^2\left<u,iv\right>-\mu \Im (uv) -\frac{1-\mu}{2}-\mu \frac{|v|^2}{|2v^2-1|}+c|v|^2,
 \end{aligned}
 \end{equation} and there is a two-to-one correspondence between a centrally symmetric sphere-like component $\widetilde\Sigma_{\mu,c}\subset K_{\mu,c}^{-1}(0)$ and $\dot \Sigma_{\mu,c}$, up to collisions.

We consider the re-scaled coordinates
\begin{equation}\label{hat}
v=\hat v\sqrt{1-\mu} \mbox{  and  } u=\hat u\sqrt{1-\mu},
\end{equation} with  Hamiltonian
\begin{equation}\label{Khat}
\begin{aligned}
&\hat  K_{\mu,c}(\hat v,\hat u) := \frac{1}{1-\mu}K_{\mu,c}(v,u)\\ & = \frac{1}{2}|\hat u|^2+2(1-\mu)|\hat v|^2\left<\hat u,i \hat v\right>-\mu \Im (\hat u \hat v)
 -\frac{1}{2}-\mu \frac{|\hat v|^2}{|2(1-\mu)\hat v^2-1|}+c|\hat v|^2.
\end{aligned}
\end{equation}
The component $\widetilde \Sigma_{\mu,c}\subset K_{\mu,c}^{-1}(0)$ gets re-scaled and we  denote it by  $\hat \Sigma_{\mu,c}\subset \hat K_{\mu,c}^{-1}(0)$.


Taking $\mu \to 1^-$ we see from \eqref{Khat} that $\hat \Sigma_{\mu,c}$ converges in the $C^\infty$ topology to a hypersurface  satisfying
\begin{equation}\label{level1}
\frac{1}{2}|\hat u|^2-\Im (\hat u \hat v) +(c-1)|\hat v|^2 = \frac{1}{2}.
\end{equation}

In order to have a better picture of the hypersurface in \eqref{level1}, we denote, for simplicity,  $\hat v=\hat v_1 + i\hat v_2$ and $\hat u=\hat u_1 + i \hat u_2$. Then \eqref{level1} is equivalent to
\begin{equation}\label{level2}
(\hat u_1-\hat v_2)^2+(\hat u_2-\hat v_1)^2 + 2\left(c-\frac{3}{2}\right)(\hat v_1^2+\hat v_2^2)= 1.
\end{equation} Taking new coordinates $(w=w_1+iw_2,z=z_1+iz_2)\in \C \times \C$
$$
\begin{aligned}
w_1=\hat u_1-\hat v_2 \mbox{  and  } z_1=\hat v_1 \sqrt{2\left(c-\frac{3}{2}\right)} ,\\
w_2=\hat u_2-\hat v_1 \mbox{  and  } z_2=\hat v_2 \sqrt{2\left(c-\frac{3}{2}\right)} ,
\end{aligned}
$$ which are symplectic up to a constant factor, we see that \eqref{level2} is equivalent to $$w_1^2+w_2^2+z_1^2+z_2^2=1.$$

We conclude that the regularized Hamiltonian flow on $\hat \Sigma_{\mu,c}$ converges smoothly to the standard Reeb flow on $(S^3,\xi_0)$ as $\mu \to 1^-$. Its orbits are the Hopf fibers. Since the projection of the retrograde orbit winds once around $0\in \C$ in $q$-coordinates, it is doubly covered by a simple closed orbit $P_{\mu,c}\subset\widetilde\Sigma_{\mu,c}$, which in $z$-coordinates winds once around $0\in \C$. Hence, $P_{\mu,c}$ converges smoothly to a Hopf fiber in $(w,z)$ and, in particular, it is unknotted and has self-linking number $-1$. The dynamical convexity of the Hamiltonian flow on $\hat \Sigma_{\mu,c}$ and Theorem~\ref{thm_hryn_dyn_convex} imply that it is the boundary of a disk-like global surface of section. In view of Theorem \ref{thm_annulus_dyn_convex}, we may assume that this global section descends to a rational disk-like global section on $\Sigma_{\mu,c} = \widetilde\Sigma_{\mu,c}/\{\pm1\}$.

\section{Transverse foliations}\label{sec_foliations}

We introduce transverse foliations adapted to a $3$-dimensional flow. This generalizes the notion of open books and global sections.

\begin{defn} Let $\phi^t$ be a smooth flow on an oriented closed $3$-manifold $M$. A transverse foliation for $\phi^t$ is formed by:
\begin{itemize}
\item[(i)] A finite set $\P$ of primitive periodic orbits of $\phi^t$, called binding orbits.
\item[(ii)] A smooth foliation of $M\setminus \cup_{P\in \P} P$ by properly embedded surfaces. Every leaf is transverse to $\phi^t$ and has an orientation induced by $\phi^t$ and $M$. For every leaf $\dot \Sigma$ there exists a compact embedded surface $\Sigma \hookrightarrow M$ so that $\dot \Sigma = \Sigma \setminus \partial \Sigma$ and $\partial \Sigma$ is a union of components of $\cup_{P\in \P} P$. An end $z$ of $\dot \Sigma$ is called a puncture. To each puncture $z$ there is an associated component $P_z\in\P$ of $\partial \Sigma$ called the asymptotic limit of $\dot \Sigma$ at $z$.
A puncture $z$ of $\dot \Sigma$ is called positive if the orientation on $P_z$ induced by $\dot \Sigma$ coincides with the orientation induced by $\phi^t$. Otherwise $z$ is called negative.
\end{itemize}
\end{defn}

The following theorem is a seminal result on the existence of transverse foliations for Reeb flows on the tight $3$-sphere. It is based on pseudo-holomorphic curve theory in symplectic cobordisms.

\begin{thm}[Hofer-Wysocki-Zehnder \cite{fols}]\label{sfef}
Let $\phi^t$ be a nondegenerate Reeb flow on $(S^3,\xi_0)$. Then $\phi^t$ admits a transverse foliation. The binding orbits have self-linking number $-1$ and their Conley-Zehnder indices are $1$, $2$ or $3$. Every leaf $\dot \Sigma$ is a punctured sphere and has precisely one positive puncture. One of the following conditions holds:
\begin{itemize}
\item The asymptotic limit of $\dot \Sigma$ at its positive puncture has Conley-Zehnder index $3$ and the asymptotic limit of $\dot \Sigma$ at any negative puncture has Conley-Zehnder index $1$ or $2$. There exists at most one negative puncture whose asymptotic limit has Conley-Zehnder index $2$.
\item The asymptotic limit of $\dot \Sigma$ at its positive puncture has Conley-Zehnder index $2$ and the asymptotic limit of $\dot \Sigma$ at any negative puncture has Conley-Zehnder index $1$.
\end{itemize}
\end{thm}

The open books with disk-like pages constructed in~\cite{char1,char2,convex,hryn,hryn2,HLS,HS1} for Reeb flows on $(S^3,\xi_0)$ are particular cases of transverse foliations with a single binding orbit. The main obstruction for the existence of such an open book with a prescribed binding orbit $P$ is the presence of closed orbits with Conley-Zehnder index $2$ which are unlinked to $P$. One particular transverse foliation of interest which deals with such situations is the so called 3-2-3 foliation.

\begin{figure}[ht] \label{fig_323}
    \centering
    \includegraphics[width=0.70\textwidth]{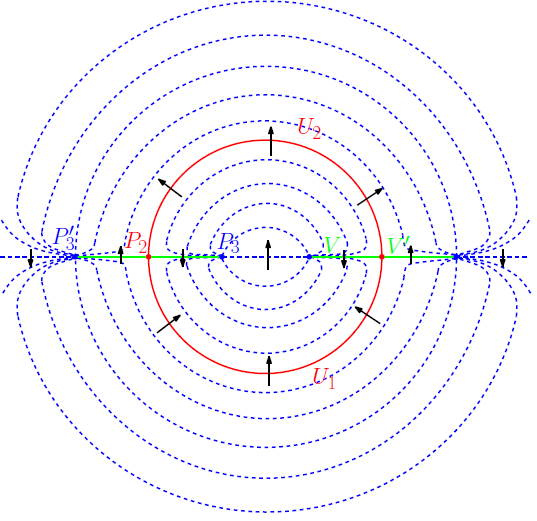}
\caption{A 3-2-3 foliation.}
\end{figure}

\begin{defn}
A 3-2-3 foliation for a Reeb flow $\phi^t$ on $(S^3,\xi_0)$ is a transverse foliation for $\phi^t$ with precisely three binding orbits $P_3$, $P_2$ and $P_3'$. They are unknotted, mutually unlinked and their Conley-Zehnder indices are $3,2$ and $3$ respectively. The leaves are punctured spheres and consist of
 \begin{itemize}
 \item A pair of planes $U_1$ and $U_2$, both asymptotic to $P_2$ at their positive punctures.
 \item A cylinder $V$ asymptotic to $P_3$ at its positive puncture and to $P_2$ at its negative puncture; a cylinder $V'$ asymptotic to $P_3'$ at its positive puncture and to $P_2$ at its negative puncture.
 \item A one parameter family of planes asymptotic to $P_3$ at their positive punctures; a one parameter family of planes asymptotic to $P_3'$ at their positive punctures.
 \end{itemize}
 See Figure \ref{fig_323}.
\end{defn}

The 3-2-3 foliations are natural objects to consider if one studies Hamiltonian dynamics near certain critical energy levels. Take a Hamiltonian  $H$ on $\R^4$ which has a Morse index $1$ critical point $p\in H^{-1}(0)$ of saddle-center type. We refer to \cite{PS1} for more precise definitions. It turns out that $p$ is a rest point of the flow and its center manifold is foliated by the so called Lyapunoff orbits $P_{2,E}\subset H^{-1}(E),E>0$ small. Each one of them is unknotted,  hyperbolic inside its energy level and has Conley-Zehnder index $2$.

\begin{figure}[ht]
    \centering
    \includegraphics[width=0.70\textwidth]{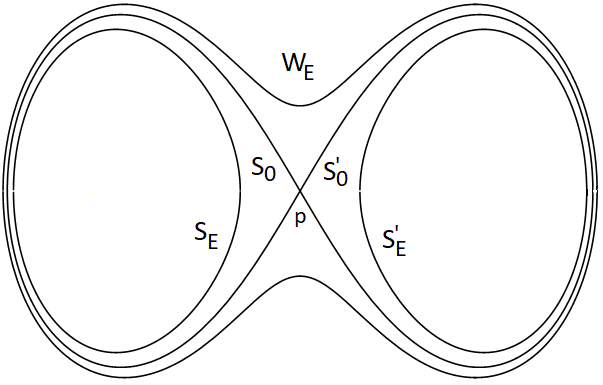}
    \caption{Strict convexity of $S_0$ and $S_0'$ implies the existence of a 3-2-3 foliation on $W_E\simeq S_E \# S_E'\simeq S^3$.}
    \label{fig_niveis}
\end{figure}

Assume that for every $E<0$ the energy level $H^{-1}(E)$ contains two sphere-like components $S_E$ and $S_E'$ which develop a common singularity at $p$ as $E \to 0^-$. This means that $S_E$ converges in the Hausdorff topology to $S_0\subset H^{-1}(0)$ as $E\to 0^-$, where $S_0$ is homeomorphic to the $3$-sphere  and contains $p$ as its unique singularity. The analog holds for $S_E'$. Therefore, $S_0 \cap S_0' =\{p\}$ and, for $E>0$ small, $H^{-1}(E)$ contains a sphere-like component $W_E$ close to $S_0\cup S_0'$. We observe that $W_E$ contains the Lyapunoff orbit $P_{2,E}$ and is in correspondence with the connected sum of  $S_E$ and $S_E'$, see Figure \ref{fig_niveis}.

\begin{defn}We say that $S_0$ is strictly convex if $S_0$ bounds a convex domain in $\R^4$ and all the sectional curvatures of $S_0\setminus \{p\}$ are positive. We say that $S_0'$ is strictly convex if analogous conditions hold.
\end{defn}

The following theorem is inspired by results in \cite{fols}.

\begin{thm}[\cite{PS1,PS2}] \label{thm_323}
If $H$ is real analytic and both $S_0$ and $S_0'$ are strictly convex then, for every $E>0$ small, the Hamiltonian flow on the sphere-like component $W_E \subset H^{-1}(E)$ admits a 3-2-3 foliation. The Lyapunoff orbit $P_{2,E}$ is one of the binding orbits and there exist infinitely many periodic orbits and infinitely many homoclinics to $P_{2,E}$ in $W_E$.
\end{thm}

One difficulty in proving Theorem~\ref{thm_323} is that there are no non-degeneracy assumptions of any kind. A criterium for checking strict convexity of the subsets $S_0$ and $S_0'$ is found in \cite{Sa1}.

The notion of 3-2-3 foliation is naturally extended to Reeb flows on connected sums $\R P^3 \# \R P^3$. In this case the binding orbits $P_3$ and $P_3'$ are non-contractible and the families of planes are asymptotic to their respective double covers. The existence of 3-2-3 foliations for Reeb flows on $\R P^3 \# \R P^3$ has interesting consequences. It is conjectured that they exist for some Hamiltonians in celestial mechanics such as Euler's problem of two centers in the plane and the planar circular restricted three body problem for energies slightly above the first Lagrange value.

A more general theory of transverse foliations for Reeb flows still needs to be developed. If one wishes to use holomorphic curves then one step is implemented by Fish and Siefring~\cite{FS}, who showed persistence under connected sums. Transverse foliations on mapping tori of disk-maps were constructed by Bramham~\cite{bramham1,bramham2} to study questions about rigidity of pseudo-rotations.

\section{A Poincar\'e-Birkhoff theorem for tight Reeb flows on $S^3$}\label{sec_PB}

Poincar\'e's last geometric theorem is nowadays known as the Poincar\'e-Birkhoff theorem. In its simplest form it is a statement about fixed points of area-preserving annulus homeomorphisms
\begin{equation*}
f:\R/\Z\times[0,1] \to \R/\Z\times[0,1]
\end{equation*}
preserving boundary components. The map $f$ can be lifted to the universal covering $\R\times[0,1]$. Let us denote projection onto the first coordinate by $p:\R\times[0,1]\to\R$. Then $f$ is said to satisfy a {\it twist condition} on the boundary if it admits a lift to the universal covering $F:\R\times[0,1]\to\R\times[0,1]$ 
such that the rotation numbers
\begin{equation*}
\lim_{n\to\infty} \frac{p\circ F^n(x,0)}{n} \qquad \lim_{n\to\infty} \frac{p\circ F^n(x,1)}{n}
\end{equation*}
differ. We call the open interval $I$ bounded by these numbers the twist interval.

\begin{thm}[Poincar\'e-Birkhoff~\cite{Bi1,Po}]
If $I\cap\Z\neq\emptyset$ then $f$ has at least two fixed points.
\end{thm}

In~\cite{franks2,franks3} John Franks has proved strong generalizations of this theorem.

Poincar\'e~\cite{Po} found annulus-like global surfaces of section for the PCR3BP for energies below the lowest critical value of the Hamiltonian, and when the mass is almost all concentrated in the primary around which the satellite moves. The boundary orbits form a Hopf link in the three-sphere (after regularization). The Poincar\'e-Birkhoff theorem applies to the associated return map and proves the existence of infinitely many periodic orbits.


One also finds such pair of orbits for the Hamiltonian flow on a smooth, compact and strictly convex energy level inside $(\R^4,\omega_0)$. In fact, the fundamental result of Hofer, Wysocki and Zehnder~\cite{convex} provides an unknotted periodic orbit $P_0$ that bounds a disk-like global surface of section. Brouwer's translation theorem yields a second periodic orbit $P_1$ simply linked to $P_0$, but much more can be said. Both orbits bind open book decompositions whose pages are disk-like global surfaces of section. It turns out that the following statement holds: the flow is equivalent to a Reeb flow on $(S^3,\xi_0)$ in such a way that $P_0\cup P_1$ corresponds to a link transversely isotopic to the standard Hopf link $$ \widetilde l_0 = \{(x_1,y_1,x_2,y_2)\in S^3 \mid x_1=y_1=0 \ \text{or} \ x_2=y_2=0\} $$ If the fixed point corresponding to $P_1$ is removed from the disk-like global section spanned by $P_0$, then we obtain a diffeomorphism of the open annulus that preserves a standard area-form and can be continuously extended to the boundary. It is interesting to study the twist condition for this map. We need to consider the transverse rotation numbers $\theta_i=\rho^{\rm Seifert}(P_i)$, $i=0,1$, with respect to Seifert surfaces (disks). In terms of Conley-Zehnder indices, these can be read as
\begin{equation}
1+\theta_0 = \lim_{n\to\infty} \frac{\CZ(P_0^n)}{2n} \qquad 1+\theta_1 = \lim_{n\to\infty} \frac{\CZ(P_1^n)}{2n}
\end{equation}
where $\CZ(P_i^n)$ denotes the Conley-Zehnder index of the $n$-iterated orbit $P_i$ computed with respect to a global trivialization of $\xi_0$. The open book singles out a lift of the map to the strip such that the rotation numbers on the boundary are precisely $1/\theta_0$ and $\theta_1$. The Poincar\'e-Birkhoff theorem proves a non-trivial statement: \\


{\it If $\theta_1 \neq 1/\theta_0$ then there are infinitely many periodic orbits in the complement of $P_0\cup P_1$. These orbits are distinguished by their homotopy classes in the complement of $P_0\cup P_1$.} \\

\begin{figure}[ht!!]
  \centering
  \includegraphics[width=0.75\textwidth]{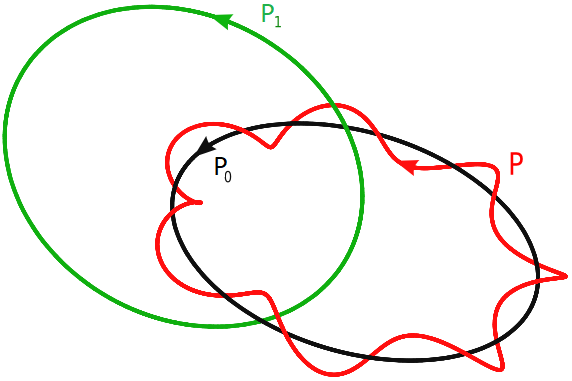}
  \caption{A Hopf link $P_0 \cup P_1$ and a closed orbit $P$ with ${\rm link}(P,P_0)=7$ and ${\rm link}(P,P_1)=1$.}
  \label{PB}
\hfill
\end{figure}


One motivation for the main result of this section is to study the possibility of extending the above result to situations where neither $P_0$ nor $P_1$ bound global sections. Before the statement we need some notation. The term {\it Hopf link} will be used for any transverse link in $(S^3,\xi_0)$ that is transversely isotopic to the standard Hopf link $\widetilde l_0$. Given non-zero vectors $u,v\in\R^2$ in the complement of the third quadrant, we write $u>v$ (or $v<u$) if the argument of $u$ is larger than that of $v$ in the counter-clockwise sense.


\begin{thm}[\cite{HMS}]\label{thm_HMS}
Consider a Reeb flow on $(S^3,\xi_0)$ that admits a pair of periodic orbits $P_0,P_1$ forming a Hopf link. Denote by $\theta_0,\theta_1$ their transverse rotation numbers computed with respect to Seifert surfaces. If $(p,q)$ is a pair of relatively prime integers satisfying $(\theta_0,1) < (p,q) < (1,\theta_1)$ or $(1,\theta_1) < (p,q) < (\theta_0,1)$
then there is a periodic orbit $P\subset S^3\setminus (P_0\cup P_1)$ such that $p={\rm link}(P,P_0)$ and $q={\rm link}(P,P_1)$.
\end{thm}

The main tools in the proof are the contact homology theory introduced by Momin~\cite{Mo} and the intersection theory of punctured holomorphic curves in dimension four developed by Siefring~\cite{siefring}. See Figure~\ref{PB}.

Another source of motivation for Theorem~\ref{thm_HMS} is a result due to Angenent~\cite{angenent} which we now recall. It concerns geodesic flows on Riemannian two-spheres. Let $g$ be a Riemannian metric on $S^2$, and let $\gamma:\R \to S^2$ be a closed geodesic of length $L$ parametrized with unit speed. In particular $\gamma(t)$ is $L$-periodic. Jacobi fields along $\gamma$ are characterized by the second order ODE
\begin{equation*}
y''(t) = -K(\gamma(t))y(t)
\end{equation*}
where $K$ denotes the Gaussian curvature. Given a (non-trivial) solution $y(t)$ we can write $y'(t)+iy(t) = r(t)e^{i\theta(t)}$ in polar coordinates. The Poincar\'e inverse rotation number of $\gamma$ is defined as
\begin{equation}
\rho(\gamma) = \frac{L}{2\pi} \lim_{t\to+\infty} \frac{\theta(t)}{t}
\end{equation}
The special case of the results from~\cite{angenent} that we would like to emphasize concerns the case when $\gamma$ is simple, that is, $\gamma|_{[0,L)}$ is injective. Denote by $n(t)$ a normal vector along $\gamma(t)$. Given relatively prime integers $p$, $q\neq0$ and $\epsilon>0$ small, a $(p,q)$-satellite about $\gamma$ is the equivalence class of the immersion $\alpha_\epsilon:\R/\Z\to S^2$
\begin{equation}
\alpha_\epsilon(t) = \exp_{\gamma(qtL)}\left( \epsilon\sin(2\pi pt) \ n(qtL) \right)
\end{equation}
where two immersions are equivalent if they are homotopic through immersions without self-tangencies and tangencies with $\gamma$.

\begin{thm}[Angenent~\cite{angenent}]\label{thm_angenent}
If a rational number $p'/q'$ strictly between $\rho(\gamma)$ and~$1$ is written in lowest terms then there exists a closed geodesic which is a $(p',q')$-satellite about $\gamma$.
\end{thm}

The above statement alone does not make justice to the non-trivial work of Angenent. He uses the curve-shortening flow to construct ``Conley pairs'' for the anti-gradient flow of the length functional associated to closed geodesics that realize a given flat knot type. The equivalence class of satellites about $\gamma$ explained above is an example of flat knot type relative to $\gamma$. His work is motivated by a question of Hofer (Oberwolfach, 1993) who asked if one could apply the Floer homology construction to curve shortening.

Let us explain the connection between theorems~\ref{thm_HMS} and~\ref{thm_angenent}. The unit tangent bundle $T^1S^2 = \{v\in TS^2 \mid g(v,v)=1 \}$ admits a contact form $\lambda_g$ whose Reeb flow coincides with the geodesic flow. It is given by the restriction to $T^1S^2$ of the pull-back of the tautological $1$-form on $T^*S^2$ by the associated Legendre transform. The $L$-periodic orbits $\dot\gamma(t)$ and $-\dot\gamma(-t)$ form a link $l_\gamma$ on $T^1S^2$ transverse to the contact structure $\ker \lambda_g$. There exists a double cover $S^3 \to T^1S^2$ that pulls back the Reeb flow of $\lambda_g$ to a Reeb flow on $(S^3,\xi_0)$. Moreover, it pulls back the link $l_\gamma$ to a Hopf link consisting of periodic orbits $P_0\cup P_1$ just like in the statement of Theorem~\ref{thm_HMS}. The main and simple observation is that $\rho(\gamma)\neq1$ forces the vectors $(\theta_0=2\rho(\gamma)-1,1)$, $(1,\theta_1=2\rho(\gamma)-1)$ to span a non-trivial sector! Then Theorem~\ref{thm_HMS} captures the contractible $(p',q')$-satellites of Theorem~\ref{thm_angenent} up to homotopy, and a refinement for Reeb flows on the standard $\R P^3$ (\cite[Theorem~1.9]{HMS}) captures all the $(p',q')$-satellites of Theorem~\ref{thm_angenent} up to homotopy. Of course, we do not hope to capture geodesics up to equivalence of satellites because Theorem~\ref{thm_HMS} deals with more general flows than those dealt by Theorem~\ref{thm_angenent}. For instance, it handles non-reversible Finsler geodesic flows with a pair of closed geodesics homotopic to a pair of embedded loops through immersions without positive tangencies. In particular, it covers reversible Finsler metrics with a simple closed geodesic.

Finally, a pair of closed Reeb orbits forming a Hopf link is not known to exist in general for a Reeb flow on $(S^3, \xi_0)$. Each of its components is unknotted, transverse to $\xi_0$ and has self-linking number $-1$; we refer to such a closed curve as a Hopf fiber. The existence of at least one closed Reeb orbit on $(S^3,\xi_0)$ which is a Hopf fiber  is proved in~\cite{unknotted}; this is a difficult result. If $P$ is a nondegenerate closed orbit which bounds a disk-like global surface of section then $P$ is a Hopf fiber and its rotation number is $>1$. Moreover, a fixed point of the first return map, assured by Brouwer's translation theorem, determines a closed orbit $P'$ which forms a Hopf link with $P$. One may ask whether every closed orbit which is a Hopf fiber and has rotation number $>1$ admits another closed orbit forming together a Hopf link. In that direction we have the following result which may be seen as a version of Brouwer's translation theorem for Reeb flows on $(S^3,\xi_0)$.

\begin{thm}[\cite{HMS2}]\label{thm_simply_linked}
Assume that a Reeb flow on $(S^3,\xi_0)$ admits a closed Reeb orbit $P$ which is a Hopf fiber. If the transverse rotation number of $P$ satisfies
$$
\rho(P) \in (1,+\infty)\setminus \left\{1+\frac{1}{k}: k\in \N\right\},
$$ then there exists a closed orbit $P'$ which is simply linked to $P$.
\end{thm}

The rotation number $\rho(P)$ is computed with respect to a global trivialization. The closed orbit $P'$ in Theorem \ref{thm_simply_linked} is not even known to be unknotted.

\appendix

\section{Infinitely many closed geodesics on a Riemannian $S^2$}

The purpose of this appendix is to describe the steps of a new proof of the existence of infinitely many closed geodesics on any Riemannian two-sphere. The argument is based on a combination of Angenent's theorem (Theorem~\ref{thm_angenent}) and the work of Hingston~\cite{hingston}, it serves as an alternative to the classical proof that combines results of Victor Bangert and John Franks. We recommend~\cite{oancea} for an account of the closed geodesic problem on Riemannian manifolds.

\begin{thm}[Bangert and Franks~\cite{bangert,franks}]\label{thm_BF}
Every Riemannian metric on $S^2$ admits infinitely many closed geodesics.
\end{thm}

We start with a remark from~\cite{hingston}. The space of embedded loops in the two-sphere carries a $3$-dimensional homology class modulo short loops. One can use Grayson's curve shortening flow to run a min-max argument over this class and obtain a special simple closed geodesic $\gamma_*$. The crucial fact here is that Grayson's curve shortening flow preserves the property of being embedded. The sum of the Morse index and the nullity of $\gamma_*$ is larger than or equal to $3$, in particular $\rho(\gamma_*)\geq1$.

If $\rho(\gamma_*)=1$ (Hingston's non-rotating case) then $\gamma_*$ is a very special critical point of the energy functional. The growth of Morse indices under iterations of $\gamma_*$ follows a specific pattern. Index plus nullity of $\gamma_*$ is equal to $3$, and if $\gamma_*$ is isolated then its local homology is non-trivial in degree $3$. If every iterate of $\gamma_*$ is isolated then the analysis of~\cite{hingston} shows that there are infinitely many closed geodesics. If some iterate $\gamma_*$ is not isolated then already there are infinitely many closed geodesics. Hence we are left with the case $\rho(\gamma_*)>1$, but this is covered by theorems~\ref{thm_HMS} or~\ref{thm_angenent} independently. Theorem~\ref{thm_BF} is proved.

\begin{rem}
It was explained in section~\ref{sec_PB} that Theorem~\ref{thm_HMS} generalizes Angenent's result to a larger class of flows. These include non-reversible Finsler metrics on~$S^2$, and one is tempted to try to push the above strategy. However, the Katok examples show that some non-reversible Finsler metrics on $S^2$ have only two closed geodesics. These metrics exist even $C^\infty$-close to the round Riemannian metric. Conjecturally there are either two or infinitely many closed geodesics.
\end{rem}

The work of Hingston~\cite{hingston,hingston2} triggered many developments. Critical points with the above properties, originally called topologically degenerate and later called symplectically degenerate maxima by Viktor Ginzburg, play key role in her proof of the Conley conjecture on standard symplectic tori~\cite{hingston_conley}: {\it Hamiltonian diffeomorphisms on standard symplectic tori must have infinitely many periodic points.} Under the generic assumption that all fixed points have no roots of unit as a Floquet multiplier the conjecture had been settled by Conley and Zehnder in~\cite{CZ_CC}, and extended to all symplectically aespherical manifolds by Salamon and Zehnder in~\cite{SZ} using Floer theory. Finally Viktor Ginzburg used Floer theory and Hinsgton's strategy to prove the following result.

\begin{thm}[Ginzburg~\cite{G}]
A Hamiltonian diffeomorphism on a symplectically aespherical closed symplectic manifold has infinitely many periodic points.
\end{thm}

Many versions of the Conley conjecture have been proved. For instance, work of Mazzucchelli~\cite{Mazz1} proves that $1$-periodic Tonelli Lagrangians over a closed base with complete Euler-Lagrange flow must have infinitely many closed orbits with integral period. See also Lu's paper~\cite{Lu}.
Work of Hein~\cite{hein} deals with the case of irrational symplectic manifolds.

Until today geometers and dynamicists struggle to understand precisely the class of phase spaces where the Conley conjecture holds. This is an active topic of research. See~\cite{GGM} for results on the analogue problem for Reeb flows.


\end{document}